%% file: linesearch_v8.tex
\documentclass{article}                                                          

\usepackage{epsfig} 
\usepackage{amsmath} 
\usepackage{amssymb}  
\usepackage{url}

\newtheorem{thm}{Theorem}

\newtheorem{lem}{Lemma}
\newtheorem{rem}{Remark}

\newenvironment{pf}{\smallbreak\noindent{\it Proof. }}{\hfill$\Box$\smallbreak}

\usepackage{tikz}
\usepackage{pgfplots}
\usetikzlibrary{calc}

\newcommand{\id}{I}
\newcommand{\reals}{\mathbb{R}}
\newcommand{\prox}{{\rm{prox}}}
\newcommand{\argmin}{\mathop{\rm argmin}}
\newcommand{\alphanom}{\bar{\alpha}}
\newcommand{\alphak}{\alpha_k}

\newcommand{\marker}[2]
{
  \pgfgettransformentries{\myxscale}{\@tempa}{\@tempa}{\myyscale}{\@tempa}{\@tempa}
  \draw[thick] ($#1+0.08*(1/\myxscale,1/\myyscale)$)--($#1-0.08*(1/\myxscale,1/\myyscale)$);
  \draw[thick] ($#1+0.08*(-1/\myxscale,1/\myyscale)$)--($#1-0.08*(-1/\myxscale,1/\myyscale)$);
}

\title{\LARGE \bf
Line Search for Averaged Operator Iteration
}

\date{}

\author{Pontus Giselsson, Mattias F{\"a}lt, and Stephen Boyd}

\begin{document}

\maketitle
\thispagestyle{empty}
\pagestyle{empty}

\begin{abstract}
Many popular first order algorithms for convex optimization, such as
forward-backward splitting, Douglas-Rachford splitting,
and the alternating direction method of multipliers (ADMM), can be 
formulated as averaged iteration of a nonexpansive mapping.
In this paper we propose a line search for averaged iteration
that preserves the theoretical convergence guarantee, while often accelerating
practical convergence.
We discuss several general cases in which the additional computational
cost of the line search is modest compared to the savings obtained.
\end{abstract}


\section{Introduction}

First-order algorithms such as forward-backward splitting, Douglas-Rachford 
splitting, and the alternating direction methods of multipliers (ADMM)
are often used for large-scale convex optimization.
While the theory tells us that these methods converge, practical convergence
can be very slow for some problem instances.
One effective method to reduce the number of iterations is to precondition the
problem data. This approach has been extensively studied in the literature and
has proven very successful in practice; see, e.g.,
\cite{Benzi_precond,Bramble_Uzawa,Hu_nonlin_Uzawa,GhadimiADMM,gisBoydAut2014metric_select,gisBoydTAC2014metric_select}
for a limited selection of such approaches.

Another general approach to improving practical efficiency is to carry 
out a line search, i.e., to first compute a tentative next iterate and then to
select the next iterate on the ray from the current iterate 
passing through the tentative iterate.
Typical line searches are based on some readily computed 
quantity such as the function value or norm of the gradient or residual.
A well designed line search preserves the theoretical convergence of the base method, 
while accelerating the practical convergence.
Line search is widely used in gradient descent or Newton methods;
see \cite{Boyd2004,Nocedal}.  These line search methods
cannot be applied to all first-order
methods mentioned above, however, since in general there is no readily computed
quantity that is decreasing.
(The convergence proofs for these methods typically rely on quantities related to the 
distance to an optimal point, which cannot be evaluated while the algorithm
is running.)
In this paper we propose a general line search
scheme that is applicable to most first-order convex optimization methods,
including those mentioned above whose convergence proofs are not
based on the decrease of an observable quantity.

We exploit the fact that many first-order optimization algorithms can be viewed
as averaged iterations of some nonexpansive operator, i.e., they can be written
in the form 
\begin{align} 
x^{k+1}=(1-\alphanom)x^k+\alphanom Sx^k = x^k + \alphanom (Sx^k-x^k),
\label{eq:averaged_iter} \end{align} 
where $\alphanom\in(0,1)$ and $S~:~\reals^n\to\reals^n$ is nonexpansive,
i.e., it satisfies $\|Su-Sv\|_2 \leq \|u-v\|_2$ for all $u,v$.
The superscript $k$ denotes iteration number.
The middle expression shows that the next point is a weighted average
of the current point $x^k$ and $Sx^k$.
The expression on the righthand side of \eqref{eq:averaged_iter} shows
that the iteration can be interpreted as a taking a step 
of length $\alphanom$
in the direction of the fixed-point residual $r^k=Sx^k-x^k$.
Assuming a fixed-point exists, the iteration \eqref{eq:averaged_iter}
converges to the set of fixed-points.

In this paper we will show how steps sometimes much larger than 
$\alphanom$ can be taken, which typically accelerates practical convergence.
This iteration has the form
\begin{align} 
x^{k+1}= x^k + \alpha_k (Sx^k-x^k),
\label{eq:averaged_iter_ls} 
\end{align} 
where $\alpha_k>0$ is chosen according to line search rules described below.
We refer to $\alpha_k$ as the \emph{step length} in the $k$th iteration,
and $\alphanom$ as the \emph{nominal step length}.
The choice $\alpha_k = \alphanom$ recovers the basic averaged iteration
\eqref{eq:averaged_iter}.
We refer to the selection of $\alpha_k$ as a line search, since we are
selecting the next iterate as a point on the line or ray
passing through $x^k$ in the direction of the residual.

The merit function used to accept
a step length $\alpha_k$ in the line search is the norm of the 
fixed-point residual $\|r\|_2=\|Sx-x\|_2$. To evaluate
this merit function for a candidate point, we must compute $Sx$, 
which corresponds to the dominant cost of taking a full iteration of 
the nominal algorithm.
In the general case, then, the line search is computationally expensive,
and there is a trade-off between the cost of the line search (which depends on
the number of candidate points examined), and the savings in iterations
due to the line search.  
But we have identified many common and interesting problem and algorithm
combinations for which the fixed-point residual can be computed at low
additional cost along the candidate ray.
In these situations, performing one iteration with line
search is roughly as expensive as performing one standard iteration of the
nominal algorithm, so the additional cost of the line search is minimal.
This happens when the nonexpansive operator $S$ can be written
as $S=S_2S_1$ where $S_1~:~\reals^n\to\reals^n$ is affine and $S_2~:~\reals^n\to\reals^n$ is relatively cheap to
evaluate.


The paper is organized as follows. 
In Section~\ref{sec:main}, we state the line search method and prove its convergence.
In Section~\ref{sec:comp_cost}, we show that the line search can be carried out
efficiently when $S=S_2S_1$ and $S_2$ is cheap to evaluate and $S_1$ is affine. 
In Section~\ref{sec:opt_algs}, we show how to implement the line search for some 
popular algorithms.
Finally, in Section~\ref{sec:num_examples} we provide numerical examples 
that show the efficiency of the proposed line search.

\section{The line search method}
\label{sec:main}

\subsection{Line search test}
The line search method first computes the nominal next iterate $\bar x^k$
according to the basic averaged iteration \eqref{eq:averaged_iter},
and then (possibly) selects a different value of $\alpha_k$.
The algorithm has the following form.
\begin{align}
\label{eq:rk}r^k&:=Sx^k-x^k\\
\label{eq:xknom}\bar x^k&:=x^k+\alphanom r^k\\
\label{eq:rknom}\bar r^k&:=S \bar x^k - \bar x^k\\
\label{eq:xk+1}x^{k+1} &:=x^k+\alphak r^k
\end{align}
In the first step we compute the current residual, in the second step
we compute the nominal next iterate, and in the third step we 
compute the nominal next residual.
In the last step, we form the actual next iterate.

In \eqref{eq:xk+1} the step length $\alphak$ must satisfy the following.
Either $\alphak = \alphanom$, i.e., we take the nominal step, or
$\alphak\in (\alphanom, \alpha^\mathrm{max}]$ is such that
\begin{align}
\|r^{k+1}\|_2 = \|Sx^{k+1}-x^{k+1}\|_2 \leq(1-\epsilon)\|\bar r^k\|_2,
\label{eq:ls_test}
\end{align}
where $\epsilon\in(0,1)$ and $\alpha^\mathrm{max} \geq \alphanom$
are fixed algorithm parameters.
Thus we either take the nominal step, or one that reduces the norm 
of the fixed point residual compared to the nominal step.

We will discuss the details of the computation and give some specific methods
to choose $\alphak$ later;
but for now we observe that to verify the line search test \eqref{eq:ls_test},
we must evaluate $r^{k+1}$, which is the first step \eqref{eq:rk} of the 
next iteration.
In a similar way, 
if we take the nominal step, i.e., choose $\alphak=\alphanom$, 
then step \eqref{eq:rknom} is the first step of the next iteration.
In either case, there is no additional computational cost.

\subsection{Convergence analysis} 

\label{sec:conv_analysis}

We analyze the proposed line search method and provide residual and iterate convergence results. All results are proven in Appendix~\ref{app:proofs}. 
\begin{thm}
Suppose that $S~:~\reals^n\to\reals^n$ is nonexpansive and 
let $\alphanom\in(0,1)$. 
Then the iteration \eqref{eq:rk}-\eqref{eq:xk+1} satisfies $\|r^k\|_2\to c$ as $k\to\infty$. 
\label{thm:norm_conv}
\end{thm}
So, the norm of the residual converges. Next, we show that the residual converges to zero if a fixed-point to $S$ exists, i.e., if ${\rm{fix}}S=\{x\in\reals^n~|~x=Sx\}\neq\emptyset$.
\begin{thm}
Suppose that $S~:~\reals^n\to\reals^n$ is nonexpansive, that ${\rm{fix}}S\neq\emptyset$, and that $\alphanom\in(0,1)$. Then the iteration \eqref{eq:rk}-\eqref{eq:xk+1} satisfies $r^k\to 0$ and $x^{k+1}\to x^k$ as $k\to\infty$. 
\label{thm:res_conv_fix_nonempty}
\end{thm}
If a fixed-point to $S$ exists, the fixed-point residual will converge to zero. Next, we establish what happens when no fixed-point to $S$ exists.
\begin{thm}
Suppose that $S~:~\reals^n\to\reals^n$ is nonexpansive, that ${\rm{fix}}S=\emptyset$, that $\inf\|Sx-x\|=c>0$, and that $\alphanom\in(0,1)$. Then the iteration \eqref{eq:rk}-\eqref{eq:xk+1} satisfies $r^k\to d$ and $x^{k+1}-x^k\to\alphanom d$ with $\|d\|=c$ as $k\to\infty$.
\label{thm:res_conv_fix_empty}
\end{thm}
This result relies heavily on \cite[Proposition~4.5]{bauschkeAffineSS} (which is a specification of more general results in \cite[Corollary~1.5]{BruckHouston1977} and \cite[Corollary 2.3]{BaillonHouston1978}). It says that, in the limit, the residual converges to a vector with smallest fixed-point residual. So the iterates converge to a line. This can, e.g., be used to devise infeasibility detection methods for these methods.

Next, we establish a rate bound for a difference of residuals.
\begin{thm}
Suppose that $S~:~\reals^n\to\reals^n$ is nonexpansive and $\alphanom\in(0,1)$. Then the iteration \eqref{eq:rk}-\eqref{eq:xk+1} satisfies
\begin{align}
  \sum_{k=0}^n\|\bar{r}^k-r^k\|_2^2\leq \frac{\alphanom}{1-\alphanom}\|r^0\|_2^2.
\label{eq:res_diff_sum}
\end{align}
Let $k_{\rm{best}}^n\in\{0,\ldots,n\}$ be the iterate $k$ (up to $n$) for which $\|\bar{r}^k-r^k\|_2$ is smallest. Then
\begin{align}
\|\bar{r}^{k_{\rm{best}}^n}-&r^{k_{\rm{best}}^n}\|_2^2
\leq \frac{\alphanom}{(n+1)(1-\alphanom)}\|r^0\|_2^2.
\label{eq:res_diff_best}
\end{align}
\label{thm:sublin_conv_rate}
\end{thm}
If $S$ is a $\delta$-contraction with $\delta\in[0,1)$, i.e., $\|Sx-Sy\|\leq\delta \|x-y\|$ for all $x,y\in\reals^n$, stronger convergence results can be obtained.
\begin{thm}
Assume that $S~:~\reals^n\to\reals^n$ is $\delta$-contractive with $\delta\in[0,1)$ and $\alphanom\in(0,1)$. Then the iteration \eqref{eq:rk}-\eqref{eq:xk+1} satisfies
\begin{align*}
\|r^{k+1}\|_2\leq (1-\alphanom+\alphanom\delta)\|r^{k}\|_2
\end{align*}
for all iterations $k$.
\label{thm:lin_conv_rate}
\end{thm}
So, the fixed-point residual converges linearly to zero (which it can since contractive operators always have a unique fixed-point).
\begin{rem}
All results in this section are stated in the Euclidean setting with the standard 2-norm. But they also hold in general finite-dimensional real Hilbert space settings.
\end{rem}

\section{Computational cost}

\label{sec:comp_cost}

The fixed-point residual must be evaluated to carry out the line search test
\eqref{eq:ls_test}. In the general case this requires us to evaluate the operator
$S$, which has the same cost as a
full iteration of the algorithm. Therefore,
in the general case it may be too expensive to evaluate many 
(or even just more than one) candidate step
lengths $\alphak$ compared to the savings in iterations due to the line search.

In this section we consider a special case in which the line search can be carried out
more efficiently, i.e.,
many candidate points along the ray can be evaluated with low additional cost. 
Suppose that $S=S_2S_1$, where $S_2~:~\reals^n\to\reals^n$ is cheap to evaluate compared to $S_1$,
and $S_1~:~\reals^n\to\reals^n$ is affine.
The algorithm \eqref{eq:rk}-\eqref{eq:xk+1} in this case becomes:
\begin{align}
\label{eq:rk_comp}r^k&:=S_2S_1x^k-x^k\\
\label{eq:xknom_comp}\bar x^k&:=x^k+\alphanom r^k\\
\label{eq:rknom_comp}\bar r^k&:=S_2S_1 \bar x^k - \bar x^k\\
\label{eq:xk+1_comp}x^{k+1} &:=x^k+\alphak r^k
\end{align}
In between \eqref{eq:rknom_comp} and \eqref{eq:xk+1_comp},
we perform the line search test \eqref{eq:ls_test},
\begin{align}
\|r^{k+1}\|_2 = \|S_2S_1x^{k+1}-x^{k+1}\|_2 \leq(1-\epsilon)\|\bar r^k\|_2,
\label{eq:ls_test_comp}
\end{align}
for multiple candidate values of $\alphak$.

We now analyze the complexity, assuming that the cost of evaluating
$S_2$, and vector-vector operations, are negligible (or at least, dominated
by the cost of evaluating $S_1$).
In one iteration with line search we need to compute $S_1x^k$ 
in \eqref{eq:rk_comp}, 
$S_1\bar x^k$ in \eqref{eq:rknom_comp}, and $S_1(x^k+\alphak r^k)$
for each candidate $\alphak$ in \eqref{eq:ls_test_comp}. 
Since $S_1$ is affine, i.e., of the form
\begin{align}
S_1(x)=Fx+h
\label{eq:S2_affine}
\end{align}
with $F\in\reals^{n\times n}$ and $h\in\reals^n$,
we have for any $\alpha$,
\begin{align*}
S_1(x^k+\alpha r^k) = F x^k + h+ \alpha Fr^k.
\end{align*}
So once we evaluate $F_2x^k$ and $F_2r^k$, we can 
evaluate $S_1(x^k + \alpha r^k)$ for any number of 
values of $\alpha$, at the cost of only vector operations.
In particular, we can evaluate $S_1 \bar x^k$ in step \eqref{eq:rknom_comp},
and $S_1 x^{k+1}$ for multiple values of $\alphak$ in the line search 
test \eqref{eq:ls_test_comp}, with no further evaluations of $S_1$.
We can express the first three steps of the algorithm as
\begin{align}
\label{eq:rk_comp2}r^k&:=S_2(Fx^k+h)-x^k\\
\label{eq:xknom_comp2}\bar x^k&:=x^k+\alphanom r^k\\
\label{eq:rknom_comp2}\bar r^k&:=S_2\left(Fx^k +h +\alphanom Fr^k\right)
- \bar x^k
\end{align}
which involves two evaluations of $F$ (and two evaluations of $S_2$), and
some vector operations. The next step is the line search, in which we evaluate the residual $r$ using
\begin{align}
r^{k+1} = S_2 \left(Fx^k+h+\alphak Fr^k\right)-(x^k+\alphak r^k)
\label{eq:res_k+1_S1_S2}
\end{align}
for $p$ candidate values of $\alphak$.  
Each of these involves a few vector operations, and one evaluation of $S_2$,
since we use the cached values of $Fr^k$ and $Fx^k$.
One iteration costs $2+p$ evaluations of $S_2$, $2$ evaluations of $F$,
and order $p$ vector operations.

Finally, as observed above, we will have 
already evaluated the step \eqref{eq:rk_comp} for the next iteration,
so one evaluation of $F$ (and $S_2$) does not count 
(or rather, counts towards the next iteration).
Thus the computational cost of one iteration with $p$ candidate 
values of $\alphak$
is one evaluation of $S_1$ (hence $F$) and $p+1$ evaluations of $S_2$.
If the cost of evaluating $S_1$ dominates the cost of 
evaluating $S_2$ (and vector operations), the computational cost
of the iteration with line search is the same as the basic 
iteration without line search.

\paragraph{A variation.} For some algorithms such as forward-backward splitting the averaged iteration \eqref{eq:averaged_iter} is more conveniently written as
\begin{align}
x^{k+1} := T_2T_1x^k
\label{eq:avg_iter_variation}
\end{align}
where $T_2~:~\reals^n\to\reals^n$ and $T_1~:~\reals^n\to\reals^n$.
So, in this case $(1-\alphanom)x^k+\alphanom S_2S_1x^k=T_2T_1x^k$. (The nominal $\alphanom$ is hidden in the composition between $T_2$ and $T_1$.)

Instead of using $S_2S_1x-x$ as residuals in \eqref{eq:rk_comp}-\eqref{eq:xk+1_comp}, we can use $\alphanom (S_2S_1x-x)=T_2T_1x-x$. An equivalent algorithm then becomes
\begin{align}
\label{eq:rk_variation}r^k&:=T_2T_1x^k-x^k\\
\label{eq:xknom_variation}\bar x^k&:=x^k+r^k\\
\label{eq:rknom_variation}\bar r^k&:=T_2T_1 \bar{x}^k- \bar{x}^k\\
\label{eq:xk+1_variation}x^{k+1}&:=x^k+\alphak r^k
\end{align}
where $\alphak\in[1,\alpha_{\max}]$.

Now, let $T_1$ be affine, i.e., of the form 
\begin{align}
T_1x=Fx+h.
\label{eq:T2_affine}
\end{align}
Then the steps \eqref{eq:rk_comp2}-\eqref{eq:rknom_comp2} (with the $x^{k+1}$ update) becomes
\begin{align}
\label{eq:rk_variation2}r^k&:=T_2(Fx^k+h)-x^k\\
\label{eq:xknom_variation2}\bar x^k&:=x^k+r^k\\
\label{eq:rknom_variation2}\bar r^k&:=T_2\left(F x^k + h +Fr^k\right)- \bar{x}^k\\
\label{eq:xk+1_variation2}x^{k+1}&:=x^k+\alphak r^k
\end{align}
The residual for the line search that is evaluated between \eqref{eq:rknom_variation2} and \eqref{eq:xk+1_variation2} is computed as
\begin{align}
r^{k+1} = T_2\left(Fx^{k}+h+\alphak Fr^k\right)-(x^k+\alphak r^k)
\label{eq:res_k+1_T1_T2}
\end{align}
for multiple candidate values of $\alphak$.

\paragraph{Evaluating affine operators.}

To evaluate the affine operator $S_1~:~\reals^n\to\reals^n$ typically involves a matrix multiplication or a matrix inversion, where the matrix is the same in all iterations. 

There are two main methods for repeated matrix inversion. 
The first is to factorize the matrix to be inverted once before the algorithm starts. Then forward and backward solves are used in every iteration. The cost of the forward and backward solves depends on the sparsity of the factors, but is typically more than $O(n)$ up to $O(n^2)$. The second option is to use an iterative method (with warm start). This requires a number of multiplications with the matrix to invert and is hence more expensive than $O(n)$. 

Assuming that the cost of evaluating $S_2~:~\reals^n\to\reals^n$ is $O(n)$, the cost of evaluating $S_1$ dominates the one of evaluating $S_2$ in this setting.

\section{Optimization algorithms}

\label{sec:opt_algs}

Many popular optimization algorithms can be implemented with the proposed line search method. In this section, we show how $S$, $S_2$ and $S_1$ (or $T_2$ and $T_1$) look for some of these. Before this, we introduce some operators.

The proximal operator associated with a proper closed and 
convex $f~:~\reals^n\to\reals\cup\{\infty\}$ is defined as
\begin{align}
\prox_{\gamma f}(z):=\argmin_{x}\{f(x)+\tfrac{1}{2\gamma}\|x-z\|_2^2\}
\label{eq:prox}
\end{align}
where $\gamma>0$. The reflected proximal operator is defined as
\begin{align}
R_{\gamma f}:=2\prox_{\gamma f}-\id.
\label{eq:rprox}
\end{align}
If $f$ is the indicator function of a nonempty closed and convex set $C$, i.e., 
\begin{align}
f(x)=\iota_{C}(x):=\begin{cases}0&{\hbox{if }} x\in C\\
  \infty & {\hbox{else}}
\end{cases}
\label{eq:indicator_fcn}
\end{align}
then the proximal operator in \eqref{eq:prox} is a projection:
\begin{align}
\prox_{\gamma f}(z)=\Pi_C(z):=\argmin_{x\in C}\|x-z\|_2
\label{eq:proj}
\end{align}
and the reflected proximal operator in \eqref{eq:rprox} is $R_{\gamma\iota_C}=R_{\iota_C}=2\Pi_C-\id$.

\subsection{Forward-backward splitting}

The forward-backward splitting method (see, e.g., \cite{CombettesFBS}) 
solves composite optimization problems of the form
\begin{align}
{\hbox{minimize }} f(x)+g(x),
\label{eq:fb_prob}
\end{align}
where $f~:~\reals^n\to\reals$ is convex and differentiable with an $L$-Lipschitz continuous
gradient $\nabla f$ and $g~:~\reals^n\to\reals\cup\{\infty\}$ is proper closed and convex.

The forward-backward algorithm for this problem is
\begin{align}
x^{k+1}:=\prox_{\gamma g}(x^k-\gamma \nabla f(x^k)),
\label{eq:fb_alg}
\end{align}
where $\gamma \in(0,\tfrac{2}{L})$ is the step size and $\prox_{\gamma g}$ is defined in \eqref{eq:prox}.

If $\gamma\in(0,\tfrac{2}{L})$, it can be shown (by combining 
\cite[Proposition~4.33]{bauschkeCVXanal}, \cite[Proposition~23.7,
Remark~4.24)(iii)]{bauschkeCVXanal}, and \cite[Proposition~2.4]{Combettes201555} or \cite[Proposition~3]{gisSIAM2015}) that 
\begin{align*}
\prox_{\gamma g}(\id-\gamma \nabla f)=(1-\alphanom)\id+\alphanom S
\end{align*}
with $\alphanom=\tfrac{2}{4-\gamma L}$, where
\[
S=(1-\tfrac{1}{\alphanom})\id+\tfrac{1}{\alphanom}\prox_{\gamma g}
(\id-\gamma\nabla f)
\]
is nonexpansive.
So, the forward-backward splitting algorithm \eqref{eq:fb_alg} is an 
averaged iteration of a nonexpansive mapping with 
$\alphanom=\tfrac{2}{4-\gamma L}$. So, if $\gamma\in(0,\tfrac{2}{L})$, we can do line search in 
forward-backward splitting.

We identify $T_2=\prox_{\gamma g}$ and $T_1=(\id-\gamma\nabla f)$ in \eqref{eq:avg_iter_variation}. With these definitions, forward-backward splitting with line search is implemented as \eqref{eq:rk_variation}-\eqref{eq:xk+1_variation}.

\paragraph{$T_1$ affine.} 

The operator $T_1=(\id-\gamma\nabla f)$ is affine if $f~:~\reals^n\to\reals$ is convex quadratic, i.e., if
\begin{align*}
f(x)=\tfrac{1}{2}x^TPx+q^Tx
\end{align*}
with $P\in\reals^{n\times n}$ positive semi-definite and $q\in\reals^n$. The operator $T_1$ becomes
\begin{align*}
T_1=(\id-\gamma P)x-\gamma q.
\end{align*}
Comparing to \eqref{eq:T2_affine}, we identify $F=\id-\gamma P$ and $h=-\gamma q$. With these $F$ and $h$, forward-backward splitting with line search can be implemented as in \eqref{eq:rk_variation2}-\eqref{eq:xk+1_variation2}. 

So a full iteration with line search needs only one multiplication with $F=(\id-\gamma P)$. If in addition $T_2=\prox_{\gamma g}$ is cheap to evaluate, one full line search iteration can be evaluated roughly at the same cost as a basic iteration of the algorithm.

\subsection{Douglas-Rachford splitting}

\label{sec:DR}

The Douglas-Rachford splitting method \cite{LionsMercier1979} solves problems of the form
\begin{align*}
{\hbox{minimize }} f(x)+g(x),
\end{align*}
where $f~:~\reals^n\to\reals\cup\{\infty\}$ and $g~:~\reals^n\to\reals\cup\{\infty\}$ are proper closed and convex. 

The algorithm is given by the following iteration
\begin{align}
\label{eq:DR1}x^{k}&:=\prox_{\gamma f}(z^k)\\
\label{eq:DR2}y^{k}&:=\prox_{\gamma g}(2x^k-z^k)\\
\label{eq:DR3}z^{k+1}&:=z^k+2\alpha(y^k-x^k)
\end{align}
where $\gamma$ is a positive scalar and $\alpha\in(0,1)$.

Using the reflected proximal operator defined in \eqref{eq:rprox}
the Douglas-Rachford algorithm can be written as
\begin{align}
z^{k+1}:=((1-\alpha)\id+\alpha R_{\gamma g}R_{\gamma f})z^k.
\label{eq:DR}
\end{align}
The reflected proximal operators $R_{\gamma g}$ and $R_{\gamma f}$ are nonexpansive \cite[Corollary~23.10]{bauschkeCVXanal}, and so is their composition $R_{\gamma g}R_{\gamma f}$.

The algorithm \eqref{eq:DR} is exactly on the form used in Section~\ref{sec:comp_cost} where $S_2=R_{\gamma g}$, $S_1=R_{\gamma f}$, $S=R_{\gamma g}R_{\gamma f}$, and $\alphanom=\alpha$. With these definitions, Douglas-Rachford with line search can be implemented as \eqref{eq:rk_comp}-\eqref{eq:xk+1_comp}.

Note that $R_{\gamma f}z^k=2x^k-z^k$ in \eqref{eq:DR1}-\eqref{eq:DR3}, $R_{\gamma g}R_{\gamma f}=2y^k-2x^k+z^k$ and the residual $r^k=R_{\gamma g}R_{\gamma f}z^k-z^k=2(y^k-x^k)$.

\paragraph{$S_1$ affine.} If $S_1=R_{\gamma f}$ is affine and $S_2=R_{\gamma g}$ is cheap to evaluate, the line search can be done almost for free, see Section~\ref{sec:comp_cost}. 

The operator $S_1=R_{\gamma f}=2\prox_{\gamma f}-\id$ is affine if $\prox_{\gamma f}$ is affine, which it is if $f$ is of the form
\begin{align*}
f(x)=\begin{cases}
\tfrac{1}{2}x^TPx+q^Tx & {\hbox{if }} Ax=b\\
\infty & {\hbox{else}}
\end{cases}
\end{align*}
with $P\in\reals^{n\times n}$ positive semi-definite, $q\in\reals^n$, $A\in\reals^{m\times n}$, and $b\in\reals^m$. (Any of the quadratic or linear functions, or the affine constraint can be removed, and the operator $S_1$ is still affine.) The proximal and reflected proximal operators of $f$ become
\begin{align*}
\prox_{\gamma f}(z)&=\begin{bmatrix}I & 0\end{bmatrix}\begin{bmatrix}
P+\gamma^{-1}I & A^T\\
A & 0
\end{bmatrix}^{-1}\begin{bmatrix}
 \gamma^{-1}z-q\\
b
\end{bmatrix}\\
R_{\gamma f}(z)&=2\prox_{\gamma f}(z)-z=2\begin{bmatrix}I & 0\end{bmatrix}\begin{bmatrix}
P+\gamma^{-1}I & A^T\\
A & 0
\end{bmatrix}^{-1}\begin{bmatrix}
 \gamma^{-1}z-q\\
b
\end{bmatrix}-z\\
&=:Fz+h
\end{align*}
where $F\in\reals^{n\times n}$ and $h\in\reals^n$.

In this situation, the first three steps of the line search algorithm are \eqref{eq:rk_comp2}-\eqref{eq:rknom_comp2} with $S_2=R_{\gamma g}$ and the residual is \eqref{eq:res_k+1_S1_S2}. As shown in Section~\ref{sec:comp_cost}, we only need one evaluation of $F$ per full iteration.

Note that in practice, the matrix $F$ is typically not stored explicitly. One alternative is to factorize $\begin{bmatrix}\begin{smallmatrix}P+\gamma^{-1} I&A^T\\A&0\end{smallmatrix}\end{bmatrix}$ before the algorithm starts. This factorization is cached and used in all consecutive iterations to compute $Fr^k$ (and $Fz^0$). Another option is to use an iterative method (with warm-start) to solve the corresponding linear system of equations.


\subsection{ADMM}

The alternating direction method of multipliers \cite{Glowinski1975,Gabay1976,BoydDistributed} solves problems of the form
\begin{align}
\begin{tabular}{ll}
minimize & $f(x)+g(z)$\\
subject to & $Ax+Bz=c$,
\end{tabular}
\label{eq:ADMM_prob}
\end{align}
where $f~:~\reals^n\to\reals\cup\{\infty\}$ and $g~:~\reals^m\to\reals\cup\{\infty\}$ are proper closed convex, and $A\in\reals^{p\times n}$, $B\in\reals^{p\times m}$, and $c\in\reals^p$. 

A standard form of ADMM (with scaled dual variable $u$ and relaxation $\alpha\in(0,1)$) is:
\begin{align}
\label{eq:ADMM1}x^{k+1} &= \argmin_{x}\{f(x)+\tfrac{\rho}{2}\|Ax+Bz^{k}-c+u^{k}\|_2^2\}\\
\label{eq:ADMM2}x_A^{k+1}&=2\alpha Ax^{k+1}-(1-2\alpha)(Bz^{k}-c)\\
\label{eq:ADMM3}z^{k+1} &= \argmin_{z}\{g(z)+\tfrac{\rho}{2}\|x_A^{k+1}+Bz-c+u^{k}\|_2^2\}\\
\label{eq:ADMM4}u^{k+1} &= u^{k}+ (x_A^{k+1}+Bz^{k+1}-c)
\end{align}
where $\alpha=\tfrac{1}{2}$ gives standard ADMM without relaxation.
This form of ADMM does not have a variable for which the algorithm is an averaged iteration of a nonexpansive mapping. 

In Appendix~\ref{app:ADMM} it is shown that ADMM is Douglas-Rachford splitting  applied to a specific problem formulation. (This is a well known fact, see, e.g., \cite{Gabay83,EcksteinPhD}.) Therefore, ADMM is $\alpha$-averaged and can be written on the form
\begin{align}
v^{k+1}=(1-\alpha)v^k+\alpha R_1R_2v^k
\label{eq:ADMM_on_DR_form}
\end{align}
where $R_1~:~\reals^p\to\reals^p$ and $R_2~:~\reals^p\to\reals^p$ are reflected proximal operators. These reflected proximal operators are given by (see \eqref{eq:rprox_p1} and \eqref{eq:rprox_p2} in Appendix~\ref{app:ADMM} where $\rho=\tfrac{1}{\gamma}$):
\begin{align}
\label{eq:R1} R_{1}(v)&=2A\argmin_{x}\{f(x)+\tfrac{\rho}{2}\|Ax-v-c\|_2^2\}-2c-v,\\
\label{eq:R2} R_{2}(v)&=-2B\argmin_{z}\{g(z)+\tfrac{\rho}{2}\|Bz+v\|_2^2\}-v.
\end{align}
The algorithm \eqref{eq:ADMM_on_DR_form} (and therefore ADMM in \eqref{eq:ADMM1}-\eqref{eq:ADMM4}) can then be implemented as (see Appendix~\ref{app:ADMM}):
\begin{align}
\label{eq:ADMM1_ls}z^k&:=\argmin_{z}\{g(z)+\tfrac{\rho}{2}\|Bz+v^k\|_2^2\}\\
\label{eq:ADMM2_ls}x^k&:=\argmin_{x}\{f(x)+\tfrac{\rho}{2}\|Ax+2Bz^k+v^k-c\|_2^2\}\\
\label{eq:ADMM3_ls}v^{k+1}&:=v^k+2\alpha(Ax^k+Bz^k-c)
\end{align}

The iteration \eqref{eq:ADMM_on_DR_form} is on the form discussed in Section~\ref{sec:comp_cost} with $S_2=R_1$, $S_1=R_2$, $S=R_1R_2$, and $\alphanom=\alpha$. With these definitions, ADMM with line search can be implemented as \eqref{eq:rk_comp}-\eqref{eq:xk+1_comp}.

Note that $R_2v^k=-2Bz^k-v^k$ in \eqref{eq:ADMM1_ls}-\eqref{eq:ADMM3_ls}, $R_1R_2v^k=2Ax^k-2c+2Bz^k+v^k$, and the residual $r^k=2(Ax^k+Bz^k-c)$ in \eqref{eq:ADMM3_ls}.

\paragraph{$R_2$ affine.} If $R_2$ is affine and $R_1$ is cheap to evaluate, then line search can be performed efficiently, see Section~\ref{sec:comp_cost}. 

The operator $R_2$ is affine if $g$ is of the form
\begin{align*}
g(z)=\begin{cases}
\tfrac{1}{2}z^TPz+q^Tz & {\hbox{if }} Lz=b\\
\infty & {\hbox{else}}
\end{cases}
\end{align*}
with $P\in\reals^{m\times m}$ positive semi-definite, $q\in\reals^m$, $A\in\reals^{s\times m}$, and $b\in\reals^s$. The operator $R_2$ in \eqref{eq:R2} becomes
\begin{align*}
R_{2}(v)&=\begin{bmatrix}-2B & 0\end{bmatrix}\begin{bmatrix}
P+\rho B^TB & L^T\\
L & 0
\end{bmatrix}^{-1}\begin{bmatrix}
 -(q+\rho B^Tv)\\
b
\end{bmatrix}-v\\
&=:Fv+h
\end{align*}
where $F\in\reals^{p\times p}$ and $h\in\reals^p$.

With these definitions of $F$ and $h$, the first three steps of ADMM with line search is \eqref{eq:rk_comp2}-\eqref{eq:rknom_comp2} with $S_2=R_1$ and the residual is \eqref{eq:res_k+1_S1_S2}. Therefore, only one application of $R_2$ (and $F$) is needed per full line search iteration, see Section~\ref{sec:comp_cost}.

Also here, the matrix $F$ is typically not stored explicitly. Instead, either a cached factorization of $\begin{bmatrix}\begin{smallmatrix}P+\rho B^TB&L^T\\L&0\end{smallmatrix}\end{bmatrix}$ or an iterative method (with warm-start) is used to compute $Fr^k$ (and $Fv^0$).

\subsection{Consensus}

The consensus algorithm \cite[Section~7]{BoydDistributed} solves problems of the form
\begin{align}
{\hbox{minimize }} f(x)=\sum_{i=1}^Nf_i(x)
\label{eq:consensus_prob}
\end{align}
where $f~:~\reals^n\to\reals\cup\{\infty\}$ and all $f_i~:~\reals^n\to\reals\cup\{\infty\}$ are proper closed and convex. An equivalent formulation is
\begin{align}
{\hbox{minimize }} f_i(x_i)+\iota_C(x_1,\ldots,x_N)
\label{eq:consensus_prob2}
\end{align}
where the consensus constraint set $C$ is 
\[
C=\{(x_1,\ldots,x_N)\in\reals^n\times\cdots\times\reals^n~|~x_1=\cdots=x_N\}
\]
and $\iota_C$ is an indicator function defined in \eqref{eq:indicator_fcn}.
That is, every $x_i\in\reals^n$ in \eqref{eq:consensus_prob2} is a local 
version of the global $x\in\reals^n$ in \eqref{eq:consensus_prob}.

We use the following formulation of the consensus algorithm:
\begin{align}
\label{eq:consensus1}x_i^{k}&:=\prox_{\gamma f_i}(2z_{\rm{av}}^k-z_i^k)\\
\label{eq:consensus2}z_i^{k+1}&:=z_i^k+(x_i^k-z_{\rm{av}}^k)
\end{align}
where $z_{\rm{av}}=\tfrac{1}{N}\sum_{i=1}^Nz_i$ is the average of the 
$z_i$'s.

This consensus algorithm is obtained by applying Douglas-Rachford splitting with $\alpha=\tfrac{1}{2}$ to \eqref{eq:consensus_prob2}. (To use ADMM as in \cite{BoydDistributed} would give an equivalent algorithm, see \cite{EcksteinPhD}, but without a variable for which the algorithm is an averaged iteration.) Therefore, it is $\tfrac{1}{2}$-averaged and can be written on the form
\begin{align*}
  {\bf{z}}^{k+1} := \tfrac{1}{2}({\bf{z}}^k+R_{\gamma f}R_{ \iota_C}{\bf{z}}^k)=\tfrac{1}{2}\left({\bf{z}}^k+R_{\gamma f}(2z_{\rm{av}}^k-{\bf{z}}^k)\right)
\end{align*}
where ${\bf{z}}=(z_1,\ldots,z_N)$.
Using local variables, it can instead be written as
\begin{align*}
z_i^{k+1} := \tfrac{1}{2}\left(z_i^k+R_{\gamma f_i}(2z_{\rm{av}}^k-z_i^k)\right)
\end{align*}
for all $i=\{1,\ldots,N\}$.

The local updates of the algorithm with line search become:
\begin{align}
\label{eq:rk_consensus}r_i^k&:= R_{\gamma f_i}(2z_{\rm{av}}-z_i^k)-z_i^k\\
  \label{eq:xknom_consensus}\bar z_i^k&:=z_i^k+ \tfrac{1}{2} r_i^k\\
\label{eq:rknom_consensus}\bar r_i^k&:=R_{\gamma f}(2\bar z_{\rm{av}}^k-\bar{z}_i^k)-\bar z_i^k\\
\label{eq:xk+1_consensus}z_i^{k+1} &:=z_i^k+\alphak r_i^k
\end{align}
where either $\alphak=\tfrac{1}{2}$, or $\alphak\in(\tfrac{1}{2},\alpha_{\max}]$ is chosen in accordance with \eqref{eq:ls_test}, i.e., such that 
\begin{align*}
\|{\bf{r}}^{k+1}\|_2\leq(1-\epsilon)\|\bar {\bf{r}}^k\|_2.
\end{align*}
where ${\bf{r}}^k=(r_1^k,\ldots,r_N^k)$.

Note that the local residual $r_i^k$ in \eqref{eq:rk_consensus} is given by $2(x_i^k-z_{\rm{av}}^k)$ in \eqref{eq:consensus2} (and similarly for $\bar{r}_i^k$ in \eqref{eq:rknom_consensus}).

The operator $R_{\iota_C}$ is always affine. Therefore, a full iteration with line search can be performed with only one evaluation of $R_{\iota_C}$, see Section~\ref{sec:comp_cost}. However, $R_{\iota_C}$ is often cheaper to evaluate than $R_{\gamma f}$. So, to evaluate a candidate point in the line search involves the costly operator $R_{\gamma f}$ and may be almost as costly as a full iteration of the algorithm.

\subsection{Alternating projection methods}
We consider the problem of finding a point in the intersection of two nonempty closed and convex sets $C$ and $D$. That is, we want to find any $x\in C\cap D$. This can equivalently be written as solving the optimization problem
\begin{align}
{\hbox{minimize }} \iota_C(x)+\iota_D(x)
\label{eq:indicator_prob}
\end{align}
where $\iota_C~:~\reals^n\to\reals\cup\{\infty\}$ and $\iota_D~:~\reals^n\to\reals\cup\{\infty\}$ are indicator functions (defined in \eqref{eq:indicator_fcn}) for $C$ and $D$ respectively .

There are numerous algorithms for finding such $x$. We focus on alternating projections and Douglas-Rachford splitting.

\paragraph{Alternating projections.}
The alternating projections \cite{vonNeumann} is given by
\begin{align}
x^{k+1}=\Pi_{C}\Pi_{D}x^k.
\label{eq:ap}
\end{align}

Since $\Pi_C$ and $\Pi_D$ are $\tfrac{1}{2}$-averaged \cite[Proposition~23.7]{bauschkeCVXanal}, the composition is $\tfrac{2}{3}$-averaged \cite[Proposition~2.4]{Combettes201555} or \cite[Proposition~3]{gisSIAM2015}. Therefore, alternating projections is an averaged iteration with $\alphanom=\tfrac{2}{3}$ and of the form $x^{k+1}=T_2T_1x^k$ where $T_2=\Pi_C$ and $T_1=\Pi_D$.

Since alternating projections is an instance of \eqref{eq:avg_iter_variation}, we can implement alternating projections with line search as \eqref{eq:rk_variation}-\eqref{eq:xk+1_variation} (with $T_2=\Pi_C$ and $T_1=\Pi_D$).





\paragraph{Douglas-Rachford.}

The problem \eqref{eq:indicator_prob} can also be solved using Douglas-Rachford splitting. The algorithm becomes
\begin{align*}
z^{k+1}=(1-\alpha)z^k+\alpha R_{\iota_C}R_{\iota_D}z^k
\end{align*}
where $\alpha\in(0,1)$. That is, we have a composition of two reflections.

This algorithm is treated in Section~\ref{sec:DR} where we identified $R_{\iota_C}=S_2$ and $R_{\iota_D}=S_1$.

\begin{rem}
Note that the $\gamma$ parameter used in standard Douglas-Rachford is not present here (since the projection is independent of this). Therefore, the only parameter to be tuned is $\alpha$, i.e., the one we perform line search over.
\end{rem}

\paragraph{$D$ affine.} 

When $D$ is affine, i.e., $D=\{x~|~Ax=b\}$, then
\begin{align*}
\Pi_{D}(x)&=\begin{bmatrix}I&0\end{bmatrix}\begin{bmatrix} I & A^T\\
A & 0\end{bmatrix}^{-1}\begin{bmatrix}x\\b\end{bmatrix},\\
R_{\iota_D}(x)&=2\Pi_D(x)-x=\begin{bmatrix}2I&0\end{bmatrix}\begin{bmatrix} I & A^T\\
A & 0\end{bmatrix}^{-1}\begin{bmatrix}x\\b\end{bmatrix}-x.
\end{align*}
Both these operators are affine. 

Assume that $\Pi_C$ (and hence $R_{\iota_C}=2\Pi_C-\id$) is cheap to evaluate. Then the line search can be implemented in alternating projections and in Douglas-Rachford splitting with almost no additional cost compared to a their basic iterations (see Section~\ref{sec:comp_cost}).

Alternating projections with line search is implemented as \eqref{eq:rk_variation2}-\eqref{eq:xk+1_variation2} with $T_2=P_C$ and $Fx+h=\Pi_D$. The residual used for the line search is \eqref{eq:res_k+1_T1_T2}. The three first steps of Douglas-Rachford with line search is \eqref{eq:rk_comp2}-\eqref{eq:rknom_comp2} with $S_2=R_{\iota_C}$ and $Fx+h=R_{\iota_D}$. The residual used for the line search is \eqref{eq:res_k+1_S1_S2}.

\subsection{Other algorithms}

There are numerous other optimization algorithms that are averaged iterations
of some nonexpansive mapping. For instance, forward-backward splitting for
solving monotone inclusion problems and for solving Fenchel dual problems, as
well as projected and standard gradient methods fit the framework. The line
search can also be used in Douglas-Rachford splitting for solving monotone
inclusion problems. Also, preconditioned ADMM methods \cite{ChambollePock} can
be interpreted as an averaged iteration of some nonexpansive mapping
\cite{HeYuan_PDconv}. The recently proposed three operator splitting method in
\cite{Davis_three_splitting} is another example. Finally, the proximal point
algorithm \cite{Rockafellar_PPA} for finding the zero of one maximally monotone
operator is an averaged iteration. Actually, an algorithm is an averaged
iteration of a nonexpansive mapping if and only if it is an instance of the proximal point method. Many of the methods
mentioned above are discussed in \cite{primer_Ryu_Boyd}.

\section{Line search variations}

There are numerous ways to create variations of the line search method. In this section, we list some that can improve practical convergence.

\paragraph{Line search activation.}

We do not need to perform line search in every iteration. Line search can be used in a subset of the iterations only. If a cheap test can indicate if a line search is beneficial, this can be used as activation rule for the line search. 

Let $v^k=x^{k}-x^{k-1}$ be the difference between consecutive iterates. We have observed that if $v^{k+1}$ and $v^k$ are almost aligned, large step lengths $\alphak$ are typically accepted. If they are not aligned, we are typically restricted to smaller $\alphak$. So, an activation rule could be that the cosine between the vectors $v^{k+1}$ and $v^k$ is large, i.e., that
\begin{align}
\frac{(v^{k+1})^Tv^k}{\|v^{k+1}\|_2\|v^k\|_2}>1-\hat{\epsilon}
\end{align}
for some small $\hat{\epsilon}>0$.

This is particularly useful for methods where the affine operator $S_1$ is not dominating (as in consensus). Even for methods where $S_1$ is dominating, this can be useful. In some cases we get fewer iterations when this activation rule is used, than if not.

\paragraph{Other candidate points.}

We are not restricted to perform the line search along the residual direction $r^k$. We can accept any candidate point $\hat{x}^{k+1}$ as the next iterate if its fixed-point residual is smaller than for the nominal point.

We introduce the residual function
\begin{align}
r(x) = Sx-x.
\label{eq:resid_fcn}
\end{align}
Then we can replace the test in \eqref{eq:ls_test} with
\begin{align}
\|r(\hat{x}^{k+1})\|_2\leq(1-\epsilon)\|r(\bar{x}^k)\|_2.
\label{eq:ls_test_new_point}
\end{align}
The full algorithm becomes
\begin{align*}
r^k&:=Sx^k-x^k\\
\bar x^k&:=x^k+\alphanom r^k\\
\bar r^k&:=S \bar x^k - \bar x^k\\
x^{k+1} &:=\begin{cases}
\hat{x}^{k+1}& {\hbox{if }} \eqref{eq:ls_test_new_point} {\hbox{ holds}}\\
x^k+\alphanom r^k & {\hbox{else}}
\end{cases}
\end{align*}
It is straightforward to verify that all convergence results for the residuals $r^k$ in Section~\ref{sec:conv_analysis} still hold in this more general setting. 

One special case is to perform line search along another direction $d^k$. Then the candidate point is $\hat{x}^{k+1}=x^k+\alphak d^k$. To evaluate the test in \eqref{eq:ls_test_new_point}, we need to compute $S_2S_1(x^{k}+\alphak d^k)$. One evaluation is in the general case as expensive as one iteration of the method. However, if $d^k=r^k$ and $S_1$ is affine, we saw in Section~\ref{sec:comp_cost} that no additional $S_1$ applications are needed to perform the line search. If the direction $d^k$ instead is a linear combination of previous residuals, i.e., $d^k=\sum_{i=0}^k\theta_ir^i$ where $\theta_i\in\reals$, also no additional applications of $S_1$ are needed due to it being affine.

\paragraph{Another line search condition.}

Here, we present another line search test that does not compare progress with a nominal step, but with the last iterate that was decided by a line search. The progress is not measured with the residual function $r$ in \eqref{eq:resid_fcn}, but with a different function $s$.

To state the line search test, we let $i_k$ be the index of the last iterate (up to the current iterate $k$) that was decided by a line search, i.e., that was not the result of a nominal step. Then any candidate point $\hat{x}^{k+1}$ can be accepted as the next iterate if the following conditions hold
\begin{align*}
\|s(\hat{x}^{k+1})\|_2\leq(1-\epsilon)\|s(x^{i_k})\|_2 {\hbox{\phantom{aaa}and\phantom{aaa}}} \|r(\hat{x}^{k+1})\|_2\leq C\|s(\hat{x}^{k+1})\|_2,
\end{align*}
where $C$ is a positive scalar, $\epsilon$ is a small positive scalar, and $r$ is the residual function in \eqref{eq:resid_fcn}. If these conditions are not satisfied, the algorithm instead takes a nominal step $x^{k+1} = x^k+\alphanom r^k$.

The convergence results in this setting become weaker. The rate results in Theorem~\ref{thm:sublin_conv_rate} and \ref{thm:lin_conv_rate} cannot be guaranteed. The results concerning the residual sequence $r^{k}$ in Theorem~\ref{thm:norm_conv}, Theorem~\ref{thm:res_conv_fix_nonempty}, and Theorem~\ref{thm:res_conv_fix_empty} can, however, be shown to hold. Let $k_0, k_1, k_2,\ldots$ be the iteration indices whose iterates have been decided by accepting a candidate line search point. Then
\begin{align*}
\|s(x^{k_p})\|_2\leq (1-\epsilon)\|s(x^{k_{p-1}})\|_2\leq(1-\epsilon)^p\|s(x^{k_0})\|_2,
\end{align*}
which implies for iteration indices $k\in[k_{p+1},k_p]$ that
\begin{align*}
\|r(x^{k})\|_2\leq \|r(x^{k_p})\|_2\leq C\|s(x^{k_p})\|_2\leq (1-\epsilon)^{p}\|s(x^{k_0})\|_2,
\end{align*}
since $\{\|r(x^k)\|_2\}$ is a nonincreasing sequence in the basic method. If the tests are satisfied an infinite number of times, then $p\to\infty$ and $\|r(x^{k})\|_2\to 0$ as $k\to\infty$. If the tests are satisfied a finite number of times (which they are if, e.g., $\inf_{x}\|Sx-x\|_2>0$), the algorithm reduces to the basic iteration after a finite number of steps. Using these insights, the proofs to the results concerning the residual $r^k$ in Theorem~\ref{thm:norm_conv}, Theorem~\ref{thm:res_conv_fix_nonempty}, and Theorem~\ref{thm:res_conv_fix_empty} can easily be modified to show that the results hold also in this setting.

To improve performance, one might want to add a condition that accepts a candidate point if there is an improvement compared to the previous iterate, i.e., if the following condition is satisfied
\begin{align*}
\|s(\hat{x}^{k+1})\|_2\leq(1-\epsilon)\|s(x^k)\|_2.
\end{align*}
This condition is, however, not needed to guarantee convergence of the method.

\section{Numerical examples}

\label{sec:num_examples}

\begin{figure}
  \centering
\input{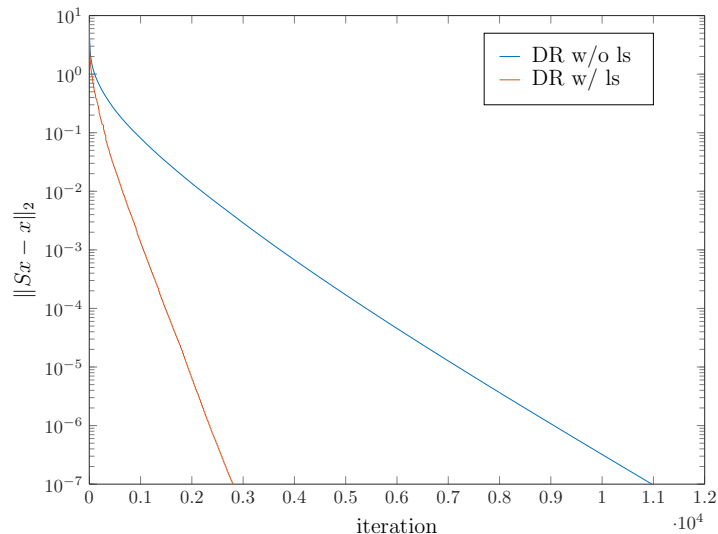}
\caption{Fixed-point residual vs iteration for Douglas-Rachford with and without line search.}
\label{fig:NNLS_rate}
\end{figure}
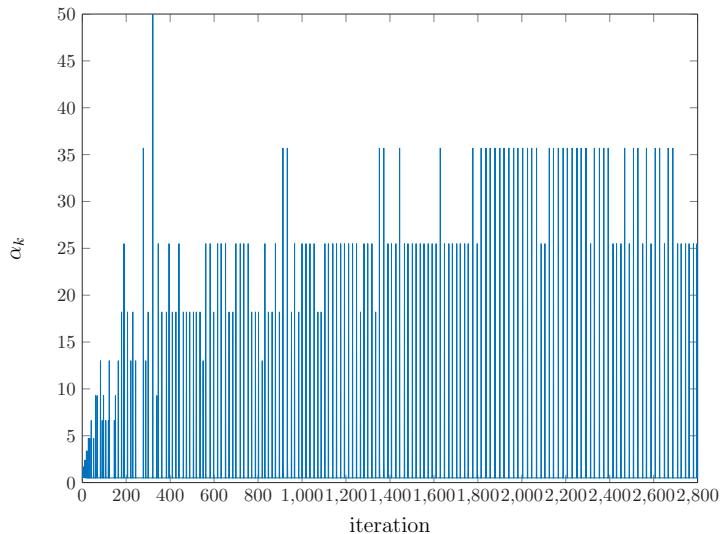
\begin{figure}
  \centering
\input{NNLS_alpha.tex}
\caption{Step length $\alphak$ vs iteration in the line search method.}
\label{fig:NNLS_alpha}
\end{figure}

\subsection{Nonnegative least squares}

To evaluate the efficiency of the line search, we solve a nonnegative least squares problem using the Douglas-Rachford algorithm. The problem is of the form
\begin{align*}
\begin{tabular}{ll}
minimize & $\|Ax-b\|_2^2$\\
subject to & $x\geq 0$
\end{tabular}
\end{align*}
where $A\in\reals^{1000\times 1000}$ is dense and $b\in\reals^{1000}$.

The entries in the data matrix $A$ are drawn from a normal distribution with zero mean and unit variance. Then, each row of $A$ is scaled with a uniformly distributed random number between 0.1 and 1.1 to worsen the conditioning of the problem. The entries in $b$ are drawn from a normal distribution with zero mean and unit variance.

To fit the Douglas-Rachford framework, we let $f(x)=\|Ax-b\|_2^2$ and $g(x)=\iota(x\geq 0)$. The operator $\prox_{\gamma f}$ is affine and the operator $\prox_{\gamma g}$ is (very) cheap to evaluate compared to $\prox_{\gamma f}$. Therefore, this problem is on the form discussed in Section~\ref{sec:comp_cost}. So an iteration with line search is just slightly more expensive than performing a basic iteration of the algorithm. 

In the line search test \eqref{eq:ls_test_comp}, we let $\epsilon=0.03$ (which may or may not be a good choice in other examples) and $\alphak$ is decided using back-tracking from $\alpha_{\max}=50$ with a factor 1/1.4 for each candidate $\alpha$. The back-tracking is stopped either when the test is satisfied, or when the candidate $\alpha\leq \alphanom$, in which case $\alphak=\alphanom$. This gives a worst case of $14$ line search test points. 

The computational cost for $\prox_{\gamma f}$ is roughly $2n^2$ after an initial matrix factorization. The cost for $\prox_{\gamma g}$ is, on the other hand, roughly $n$. To evaluate the line search test, no additional $\prox_{\gamma f}$ computations are needed. But about 10 vector additions or multiplications with scalars and one $\prox_{\gamma g}$ is needed for every candidate point (the same as in the standard algorithm). So, evaluating one candidate point costs approximately $10n$. A worst case of 14 candidate points costs $140n$ for a full line search. Comparing this to the cost for one basic iteration, $2n^2+10n$, gives, when $n=1000$, that one iteration with line search costs, in the worst case, 1.07 times a basic iteration.

Figure~\ref{fig:NNLS_rate} shows the fixed-point residual vs iteration number for Douglas-Rachford with and without line search (the Douglas-Rachford parameters are chosen to be $\alphanom=\tfrac{1}{2}$ and $\gamma= 3$). For this example, the number of iterations is reduced by roughly a factor four. The improvement in execution time is roughly the same because of the modest 7$\%$ increase in computational cost due to the line search. 

Figure~\ref{fig:NNLS_alpha} shows what values $\alphak$ that are chosen in the line search. An $\alphak=\alphanom$ corresponds to a standard Douglas-Rachford iteration. In 175 out of the 2800 iterations, an $\alphak>\alphanom$ was selected. Among these 158 had $\alphak>5$.

\def\thetainit{350}

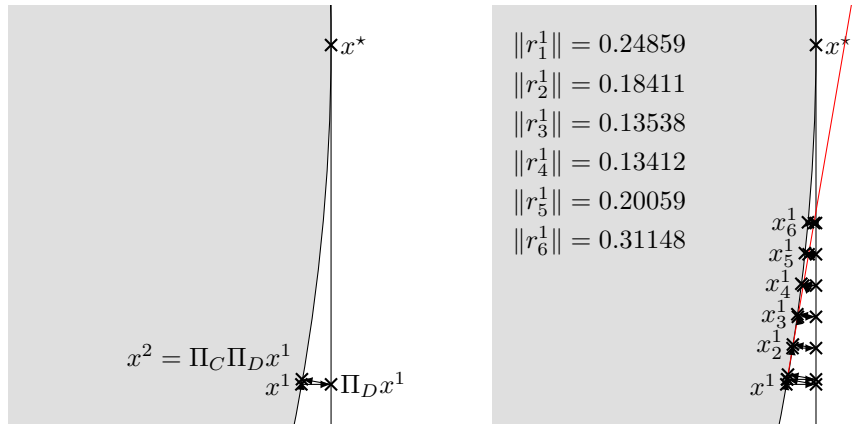
\begin{figure}
\centering
\begin{tikzpicture}


\begin{scope}[scale=26]
\clip ({cos(\thetainit)-0.15},{sin(\thetainit)-0.02}) rectangle (1.05,0.02);

\def\xline{1};
\draw[fill=gray!25!white] (0,0) circle (1);
\draw (\xline,-2)--(\xline,2);
\marker{(1,0)}{0.01}
\node[right] at (1,0) {$x^\star$};
\edef\xoneinit{cos(\thetainit)};
\edef\xtwoinit{sin(\thetainit)};

\node[above] at (0.925,{sin(\thetainit)-0.1}) {$C$};
\node[above] at (0.99,{sin(\thetainit)-0.1}) {$D$};

\pgfmathparse{1*\xoneinit}
\edef\xone{\pgfmathresult}
\pgfmathparse{1*\xtwoinit}
\edef\xtwo{\pgfmathresult}

\marker{(\xone,\xtwo)}{0.01};
\node[left] at (\xone,\xtwo) {$x^{1}$};

\edef\pxone{\xline};
\pgfmathparse{1*\xtwo}
\edef\pxtwo{\pgfmathresult}
\marker{(\pxone,\pxtwo)}{0.01};
\node[right] at (\pxone,\pxtwo) {$\Pi_Dx^{1}$};

\pgfmathparse{sqrt(\pxone*\pxone+\pxtwo*\pxtwo)};
\edef\pxnorm{\pgfmathresult};
\pgfmathparse{\pxone/\pxnorm};
\xdef\ppxone{\pgfmathresult};
\pgfmathparse{\pxtwo/\pxnorm};
\xdef\ppxtwo{\pgfmathresult};
\marker{(\ppxone,\ppxtwo)}{0.01};
\node[above left] at (\ppxone,\ppxtwo) {$x^{2}=\Pi_C\Pi_Dx^1$};

\draw[-latex] (\xone,\xtwo)--(\pxone,\pxtwo);
\draw[-latex] (\pxone,\pxtwo)--(\ppxone,\ppxtwo);
\draw[-latex] (\xone,\xtwo)--(\ppxone,\ppxtwo);

\pgfmathparse{1*\ppxone}
\xdef\xoneold{\pgfmathresult};
\pgfmathparse{1*\ppxtwo}
\xdef\xtwoold{\pgfmathresult};

\end{scope}

\end{tikzpicture}
\hspace{5mm}\newcommand\optresnorm{1000}
\newcommand\optt{1}
\begin{tikzpicture}
\input{ap_step_w_ls.tex}
\end{tikzpicture}

\caption{The left figure shows one iteration of alternating projections. The residual in this figure is $r^1=x^2-x^1$. In the right figure, an alternating projections step with line search is performed. The residual direction is shown in red. We evaluate six candidate points $x_{i}^1$, $i\in\{1,\ldots,6\}$, along this line. (The points themselves, $\Pi_D x_{i}^1$ and $\Pi_C\Pi_D x_{i}^1$ are marked with crosses in the figure.) The norm of each residual $r_{i}^1=\Pi_C\Pi_D x_{i}^1-x_{i}^1$ is printed in the figure. The 4th point $x_4^1$ has smallest residual norm. This corresponds to $\alphak=19.75$. Another option is to choose the farthest candidate point with residual norm smaller than $\|r_1^1\|$. This holds for the fifth point with $\alphak=26$. Both these choices are convergent. In this case we get closer to the intersection point by choosing the farthest point.}
\label{fig:ap_w_wo_ls}
\end{figure}

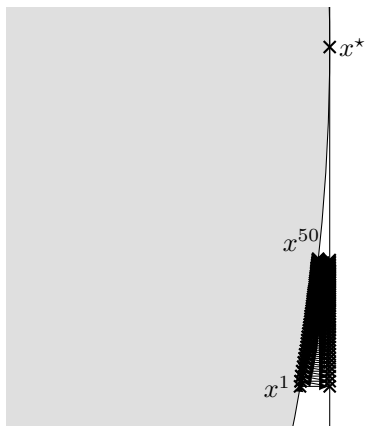
\begin{figure}
\centering
\begin{tikzpicture}


\begin{scope}[scale=26]
\clip ({cos(\thetainit)-0.15},{sin(\thetainit)-0.02}) rectangle (1.05,0.02);

\def\xline{1};
\draw[fill=gray!25!white] (0,0) circle (1);
\draw (\xline,-2)--(\xline,2);
\marker{(1,0)}{0.01}
\node[right] at (1,0) {$x^\star$};
\edef\xoneinit{cos(\thetainit)};
\edef\xtwoinit{sin(\thetainit)};
\xdef\one{1}
\node[above] at (0.925,{sin(\thetainit)-0.1}) {$C$};
\node[above] at (0.99,{sin(\thetainit)-0.1}) {$D$};
\xdef\maxiters{49}

\foreach \i in {1,2,...,\maxiters} {

\pgfmathparse{\one}
\ifx\i\pgfmathresult
\pgfmathparse{1*\xoneinit}
\edef\xone{\pgfmathresult}
\pgfmathparse{1*\xtwoinit}
\edef\xtwo{\pgfmathresult}
\else
\pgfmathparse{1*\xoneold}
\edef\xone{\pgfmathresult};
\pgfmathparse{1*\xtwoold};
\edef\xtwo{\pgfmathresult};
\fi

\marker{(\xone,\xtwo)}{0.01};

\pgfmathparse{\one}
\ifx\i\pgfmathresult
\node[left] at (\xone,\xtwo) {$x^{1}$};
\else
\fi

\edef\pxone{\xline};
\pgfmathparse{1*\xtwo}
\edef\pxtwo{\pgfmathresult}
\marker{(\pxone,\pxtwo)}{0.01};

\pgfmathparse{sqrt(\pxone*\pxone+\pxtwo*\pxtwo)};
\edef\pxnorm{\pgfmathresult};
\pgfmathparse{\pxone/\pxnorm};
\xdef\ppxone{\pgfmathresult};
\pgfmathparse{\pxtwo/\pxnorm};
\xdef\ppxtwo{\pgfmathresult};
\marker{(\ppxone,\ppxtwo)}{0.01};

\draw[-latex] (\xone,\xtwo)--(\pxone,\pxtwo);
\draw[-latex] (\pxone,\pxtwo)--(\ppxone,\ppxtwo);
\draw[-latex] (\xone,\xtwo)--(\ppxone,\ppxtwo);

\pgfmathparse{\maxiters}
\ifx\i\pgfmathresult
\pgfmathparse{\maxiters+1}
\edef\laststep{\pgfmathresult}
\node[above left] at (\pxone,\pxtwo) {$x^{\pgfmathprintnumber{\laststep}}$};
\else
\fi

\pgfmathparse{1*\ppxone}
\xdef\xoneold{\pgfmathresult};
\pgfmathparse{1*\ppxtwo}
\xdef\xtwoold{\pgfmathresult};

}

\end{scope}
\end{tikzpicture}

\caption{This figure shows 50 iterations of alternating projections. Comparing to Figure~\ref{fig:ap_w_wo_ls} reveals that roughly 50 steps of alternating projections give the same progression as one step with line search (when the farthest acceptable point is chosen) in this example. }
\label{fig:ap_many_steps}
\end{figure}

\subsection{An alternating projections example}

To visualize the line search, we solve a two dimensional feasibility problem using alternating projections. 

We want to find a point in the intersection between two sets $C=\{x\in\reals^2~|~\|x\|\leq 1\}$ and $D=\{x\in\reals^2~|~x=(x_1,x_2),x_1=1\}$. So $C$ is the unit circle, and $D$ is a vertical line that touches the boundary of $C$ at $(1,0)$. The unique intersection point is $x^\star=(1,0)$.

In Figure~\ref{fig:ap_w_wo_ls} we show one iteration of the standard alternating projections algorithm and one iteration with line search. In Figure~\ref{fig:ap_many_steps} we show 50 steps of alternating projections.

We see that the progression in 50 steps of alternating projections is roughly the same as the progression of one step with line search (when the farthest acceptable candidate point is chosen). The line search scheme does, on the other hand, compute six candidate points to advance this far. (Or really five, since the first is the basic next step.) So, we gain roughly a factor 10 in this step.

This is just a simple example where both projections are very cheap. If the cost of projecting onto the subspace is dominating the other cost of the other projection. Then the cost of performing one iteration with line search is roughly the same as the cost of one basic iteration. In such cases, we can gain a lot by performing line search.








\section{Acknowledgments}

The first author is financially supported by the Swedish Foundation for Strategic Research. The two first authors are members of the LCCC Linneaus Center at Lund University.

\bibliographystyle{plain}

\bibliography{/local/home/pontusg/Research/MPC/papers/references/references.bib}

\appendix

\section{Proofs to results in Section~\ref{sec:main}}

\label{app:proofs}

\subsection{Proof to Theorem~\ref{thm:norm_conv}}

First, we show that $\|r^k\|_2=\|x^k-Sx^k\|_2\to c$ as $k\to\infty$. We show this by considering the cases $\alphak=1$ and $\alphak>1$ separately. 

First, we consider the case $\alphak=1$. For convenience, we introduce the operator $T=(1-\alphanom)\id+\alphanom S$. Then the update for $\bar{x}^k$ in \eqref{eq:xknom} can be written as
\begin{align*}
\bar{x}^k = x^k+\alphanom(Sx^k-x^k)=(1-\alphanom)x^k+\alphanom Sx^k=Tx^k.
\end{align*} 
Noting that $\|x-Tx\|_2=\|x-(1-\alphanom)x-\alphanom Sx\|_2=\alphanom\|x-Sx\|_2$ implies
\begin{align*}
\|r^{k+1}\|_2=\|\bar{r}^k\|_2=\|\bar{x}^{k}-S\bar{x}^{k}\|_2=\tfrac{1}{\alphanom}\|\bar{x}^{k}-T\bar{x}^{k}\|_2=\tfrac{1}{\alphanom}\|Tx^k-TTx^k\|_2.
\end{align*}
Therefore, since $T$ is nonexpansive:
\begin{align}
\label{eq:no_ls_decrease}\|r^{k+1}\|_2&\leq\tfrac{1}{\alphanom}\|x^k-Tx^k\|_2=\|x^k-Sx^k\|_2=\|r^k\|_2.
\end{align}

Next, we
consider the case where $\alphak>1$. Since $\|\bar{r}^{k}\|_2\leq \|r^k\|_2$, we get from the line search test \eqref{eq:ls_test} that
\begin{align}
\label{eq:ls_decrease}\|r^{k+1}\|_2\leq (1-\epsilon)\|\bar{r}^{k}\|_2\leq (1-\epsilon)\|r^k\|_2.
\end{align}
So $\{\|r^{k}\|_2\}_{k=1}^{\infty}$ is a decreasing sequence which is bounded below (by 0). Hence it converges. This completes the proof.

\subsection{Proof to Theorem~\ref{thm:res_conv_fix_nonempty}}

Combining \eqref{eq:no_ls_decrease} and \eqref{eq:ls_decrease}, we get
\begin{align}
\label{eq:res_decrease}\|r^{k+1}\|_2\leq (1-\epsilon)^{k_0}\|r^0\|_2
\end{align}
where $k_0$ is the number of times that $\alphak$ satisfies $\alphak>1$. If $k_0\to\infty$ as $k\to\infty$, then $\|r^{k+1}\|_2\to 0$ as $k\to\infty$. On the other hand, if $k_0$ stays finite as $k\to\infty$, there exists a finite $k_{\max}$ after which the line search is not activated again. Then for $k\geq k_{\max}$, the algorithm reduces to $x^{k+1}=Tx^k$, which satisfies $\|r^k\|_2=\|x^{k}-Sx^{k}\|_2=\tfrac{1}{\alphanom}\|x^{k}-Tx^{k}\|_2\to 0$ as $k\to\infty$, see \cite[Theorem~5.14]{bauschkeCVXanal}.
This concludes the proof.

\subsection{Proof to Theorem~\ref{thm:res_conv_fix_empty}}

Combining \eqref{eq:no_ls_decrease} and \eqref{eq:ls_decrease}, we get
\begin{align}
\label{eq:res_decrease}\|r^{k+1}\|_2\leq (1-\epsilon)^{k_0}\|r^0\|_2
\end{align}
where $k_0$ is the number of times that $\alphak$ satisfies $\alphak>1$. If $k_0\to\infty$ as $k\to\infty$, then $\|r^{k+1}\|_2\to 0$ as $k\to\infty$. This is a contradiction to that $\inf\|Sx-x\|_2>0$. Hence $k_0$ must be finite and there exists a $k_{\max}$ after which the algorithm reduces to the basic averaged iteration.

Let $T=(1-\alphanom)\id+\alphanom S$, $x^{k_{\max}}=\tilde{x}_0$ and $\Delta k = k-k_{\max}$. Then a straightforward generalization of \cite[Proposition~4.5]{bauschkeAffineSS} to allow for averaged operators (instead of only firmly nonexpansive or $\tfrac{1}{2}$-averaged) gives that
\begin{align*}
\|\alphanom r^{k}-v\|=\|x^k-x^{k+1}-v\|=\|T^{\Delta k}\tilde{x}_0-T^{\Delta k+1}\tilde{x}_0-v\|\to 0
\end{align*}
for a specific $v$. Therefore $r^k\to \tfrac{1}{\alphanom}v=:d$ as $k\to\infty$. Further, $x^{k+1}-x^k=\alphanom r^k\to\alphanom d$ as $k\to\infty$. 

The $v$ is the {\emph{infimal displacement vector}} (see \cite[Fact~2.2]{bauschkeAffineSS}) that satisfies $v\in\overline{{\rm{ran}}}(\id-T)$ (i.e., $v$ is in the closure of the range of $\id-T$) and $\|v\|_2=\inf_x\|x-Tx\|_2$. Therefore $\|d\|_2=\tfrac{1}{\alphanom}\|v\|_2=\tfrac{1}{\alphanom}\inf_x\|x-Tx\|_2=\inf_x\|x-Sx\|_2$. This concludes the proof.

\subsection{Proof to Theorem~\ref{thm:sublin_conv_rate}}

We need the following lemma for this proof.
\begin{lem}
Suppose that $S~:~\reals^n\to\reals^n$ nonexpansive and that $\alphanom\in(0,1)$. Then every iteration of \eqref{eq:rk}-\eqref{eq:xk+1} satisfies
\begin{align}
\label{eq:telescope_eq}\alphanom(1-\alphanom)\|\bar{r}^k-&r^k\|_2^2
\leq \|x^k-\bar{x}^{k}\|_2^2-\|x^{k+1}-\bar{x}^{k+1}\|_2^2.
\end{align}
\label{lem:telescope_eq}
\end{lem}
\begin{pf}
Let $T=(1-\alphanom)\id+\alphanom S$. Then $T$ is $\alphanom$-averaged, and it satisfies \cite[Proposition~4.25(iii)]{bauschkeCVXanal}
\begin{align*}
\tfrac{1-\alphanom}{\alphanom}\|(\id-T)\bar{x}^k-(\id-T)x^k\|_2^2
\leq \|x^k-\bar{x}^k\|_2^2-\|Tx^k-T\bar{x}^k\|_2^2.
\end{align*}
Now, since $(\id-T)x=(\id-(1-\alphanom)\id-\alphanom S)x=\alphanom (\id-S)x$, we have $(\id-T)x^k=\alphanom r^k$ and $(\id-T)\bar{x}^k=\alphanom \bar{r}^k$. Therefore
\begin{align*}
\alphanom(1-\alphanom)\|\bar{r}^k-r^k\|_2^2
\leq \|x^k-\bar{x}^k\|_2^2-\|Tx^k-T\bar{x}^k\|_2^2.
\end{align*}
The algorithm chooses either $\alphak=\alphanom$ or $\alphak>\alphanom$.
If $\alphak=\alphanom$, we have $Tx^k=\bar{x}^k=x^{k+1}$ and $T\bar{x}^k=Tx^{k+1}=\bar{x}^{k+1}$. Therefore
\begin{align*}
  \alphanom(1-\alphanom)\|\bar{r}^k-r^k\|_2^2
&\leq \|x^k-\bar{x}^{k}\|_2^2-\|Tx^k-T\bar{x}^{k}\|_2^2\\
&= \|x^k-\bar{x}^{k}\|_2^2-\|x^{k+1}-\bar{x}^{k+1}\|_2^2.
\end{align*}
If instead $\alphak>\alphanom$, we get
\begin{align*}
\alphanom(1-\alphanom)\|\bar{r}^k-r^k\|_2^2
&\leq \|x^k-\bar{x}^k\|_2^2-\|Tx^k-T\bar{x}^k\|_2^2\\
&= \|x^k-\bar{x}^k\|_2^2-\|\bar{x}^k-T\bar{x}^k\|_2^2\\
&\leq \|x^k-\bar{x}^k\|_2^2-\tfrac{1}{(1-\epsilon)^2}\|x^{k+1}-Tx^{k+1}\|_2^2\\
&\leq \|x^k-\bar{x}^k\|_2^2-\|x^{k+1}-Tx^{k+1}\|_2^2\\
&= \|x^k-\bar{x}^k\|_2^2-\|x^{k+1}-\bar{x}^{k+1}\|_2^2
\end{align*}
where the second inequality holds due to the line search test in \eqref{eq:ls_test} and the third inequality holds since $\epsilon\in(0,1)$. Therefore \eqref{eq:telescope_eq} holds for all $k$ and the proof is complete.
\end{pf}
Now we are ready to prove the result.
A telescope summation of \eqref{eq:telescope_eq} gives
\begin{align*}
  \alphanom(1-\alphanom)\sum_{k=0}^n\|\bar{r}^k-r^k\|_2^2\leq \|x^0-\bar{x}^0\|_2^2=\alphanom^2\|r^0\|_2^2.
\end{align*}
This proves \eqref{eq:res_diff_sum}. To prove \eqref{eq:res_diff_best}, we note that $k_{\rm{best}}^n\in\{0,\ldots,n\}$ is the iteration $k$ (up till $n$) with smallest $\|\bar{r}^k-r^k\|_2$. Therefore
\begin{align*}
(n+1)\|\bar{r}^{k_{\rm{best}}^n}-r^{k_{\rm{best}}^n}\|_2^2
  \leq \sum_{k=0}^n\|\bar{r}^k-r^k\|_2^2
\leq \frac{\alphanom}{1-\alphanom}\|r^0\|_2^2.
\end{align*}
This concludes the proof.

\subsection{Proof to Theorem~\ref{thm:lin_conv_rate}}

First, we introduce $T=(1-\alphanom)\id+\alphanom S$ which is $\alphanom$-averaged, and satisfies $\|x-Sx\|_2=\tfrac{1}{\alphanom}\|x-Tx\|_2$. 
Let's consider the case when $\alphak=\alphanom$. Then $\bar{x}^k=Tx^k$ and
\begin{align*}
 \|r^{k+1}\|_2&=\|\bar{r}^k\|_2=\|\bar{x}^{k}-S\bar{x}^{k}\|_2=\tfrac{1}{\alphanom}\|\bar{x}^{k}-T\bar{x}^{k}\|_2
=\tfrac{1}{\alphanom}\|Tx^k-TTx^k\|_2\\
&=\tfrac{1}{\alphanom}\|(1-\alphanom)(x^k-Tx^k)+\alphanom (Sx^k-STx^k)\|_2.
\end{align*}
The triangle inequality gives that
\begin{align*}
 \|r^{k+1}\|_2&\leq\tfrac{1}{\alphanom}((1-\alphanom)\|x^k-Tx^k\|_2+\alphanom\|Sx^k-STx^k\|_2)\\
&\leq\tfrac{1}{\alphanom}((1-\alphanom)\|x^k-Tx^k\|_2+\alphanom\delta\|x^k-Tx^k\|_2)\\
&=\tfrac{1}{\alphanom}(1-\alphanom+\alphanom\delta)\|x^k-Tx^k\|_2\\
&=(1-\alphanom+\alphanom\delta)\|x^k-Sx^k\|_2\\
&=(1-\alphanom+\alphanom\delta)\|r^k\|_2.
\end{align*}
Next, we consider the case when $\alphak>\alphanom$. Since $\|\bar{r}^k\|_2\leq(1-\alphanom+\alphanom\delta)\|r^k\|_2$ the line search test \eqref{eq:ls_test} implies that
\begin{align*}
\|r^{k+1}\|_2\leq (1-\epsilon)\|\bar{r}^k\|_2\leq (1-\epsilon)(1-\alphanom+\alphanom\delta)\|r^k\|_2\leq (1-\alphanom+\alphanom\delta)\|r^k\|_2.
\end{align*}
That is, the algorithm is linearly convergent with factor (at most) $(1-\alphanom+\alphanom\delta)$ in both situations. This concludes the proof.

\section{ADMM derivation}

\label{app:ADMM}

In this section, we show the equivalence between the standard ADMM formulation \eqref{eq:ADMM1}-\eqref{eq:ADMM4} and the ADMM version used for line search \eqref{eq:ADMM1_ls}-\eqref{eq:ADMM3_ls}. We also show that the version used for line search, \eqref{eq:ADMM1_ls}-\eqref{eq:ADMM3_ls}, is an $\alpha$-averaged iteration of a nonexpansive mapping.

We do this by showing that the ADMM iterations can be derived by applying Douglas-Rachford splitting to a specific problem formulation. This derivation is not new \cite{Gabay83,EcksteinPhD}, but we include it here for completeness and to explicitly arrive that the ADMM variation \eqref{eq:ADMM1_ls}-\eqref{eq:ADMM3_ls} that we need for the line search.

ADMM solves problems of the form
\begin{align}
\begin{tabular}{ll}
minimize & $f(x)+g(z)$\\
subject to & $Ax+Bz=c$
\end{tabular}
\label{eq:ADMM_prob_app}
\end{align}
where $f~:~\reals^n\to\reals\cup\{\infty\}$ and $g~:~\reals^m\to\reals\cup\{\infty\}$ are proper closed convex, $A\in\reals^{p\times n}$, $B\in\reals^{p\times m}$, and $c\in\reals^{p}$.

Using image functions (that are also called infimal postcompositions) defined as
\begin{align*}
(L\triangleright \psi)(y)=\inf\{\psi(x)~|~Lx=y\}
\end{align*} 
where $L\in\reals^{n\times m}$ is a linear operator and $\psi~:~\reals^n\to\reals\cup\{\infty\}$ is a proper function, 
it is straightforward to verify that \eqref{eq:ADMM_prob_app} is equivalent to
\begin{align*}
{\hbox{minimize }} (-A\triangleright f)(-u-c)+(-B\triangleright g)(u).
\end{align*}
Let $p_1(u)=(-A\triangleright f)(-u-c)$ and $p_2(u)=(-B\triangleright g)(u)$ to get the equivalent problem
\begin{align}
{\hbox{minimize }} p_1(u)+p_2(u).
\label{eq:DR_imgfcn_prob}
\end{align}
To arrive at the standard ADMM iterations, we apply Douglas-Rachford splitting to \eqref{eq:DR_imgfcn_prob}. The algorithm becomes
\begin{align}
v^{k+1}=(1-\alpha)v^k+\alpha R_{\gamma p_1}R_{\gamma p_2}v^k
\label{eq:ADMM_avg_iter}
\end{align}
where the reflected proximal operators $R_{\gamma p_1}$ and $R_{\gamma p_2}$ are given by $R_{\gamma p_1}=2\prox_{\gamma p_1}-\id$ and $R_{\gamma p_2}=2\prox_{\gamma p_2}-\id$.
Under the assumption that the infimum over $x$ is attained in the following prox evaluation, we have
\begin{align}
\nonumber\prox_{\gamma p_1}(v)&=\argmin_{u}\{\inf_{x}\{f(x)~|~-Ax=-u-c\}+\tfrac{1}{2\gamma}\|u-v\|_2^2\}\\
\label{eq:prox_p1}&=A\argmin_{x}\{f(x)+\tfrac{1}{2\gamma}\|Ax-v-c\|_2^2\}-c.
\end{align}
The reflected proximal operator becomes
\begin{align}
R_{\gamma p_1}(v)=2A\argmin_{x}\{f(x)+\tfrac{1}{2\gamma}\|Ax-v-c\|_2^2\}-2c-v.
\label{eq:rprox_p1}
\end{align}
Again, assuming that the following infimum is attained, we get
\begin{align}
\nonumber\prox_{\gamma p_2}(v)&=\argmin_{u}\{\inf_{z}\{g(z)~|~-Bz=u\}+\tfrac{1}{2\gamma}\|u-v\|_2^2\}\\
\label{eq:prox_p2}&=-B\argmin_{z}\{g(z)+\tfrac{1}{2\gamma}\|Bz+v\|_2^2\}
\end{align}
and reflected proximal operator
\begin{align}
R_{\gamma p_2}(v)=-2B\argmin_{z}\{g(z)+\tfrac{1}{2\gamma}\|Bz+v\|_2^2\}-v.
\label{eq:rprox_p2}
\end{align}
Using the prox expressions \eqref{eq:prox_p1} and \eqref{eq:prox_p2}, and defining $\rho=\tfrac{1}{\gamma}$, we find that the Douglas-Rachford algorithm \eqref{eq:DR1}-\eqref{eq:DR3} applied to \eqref{eq:DR_imgfcn_prob} becomes
\begin{align}
\label{eq:ADMM1_ls_app}z^k&=\argmin_{z}\{g(z)+\tfrac{\rho}{2}\|Bz+v^k\|_2^2\}\\
\label{eq:ADMM2_ls_app}x^k&=\argmin_{x}\{f(x)+\tfrac{\rho}{2}\|Ax+2Bz^k+v^k-c\|_2^2\}\\
\label{eq:ADMM3_ls_app}v^{k+1}&=v^k+2\alpha (Ax^k+Bz^k-c)
\end{align}
This is exactly the iteration \eqref{eq:ADMM1_ls}-\eqref{eq:ADMM3_ls} which is used in the line search. This algorithm is equivalent to ADMM, but keeps the $v^k$ variables in which the algorithm can be interpreted as an averaged iteration of a nonexpansive mapping, see \eqref{eq:ADMM_avg_iter}.

To derive the ADMM iterations \eqref{eq:ADMM1}-\eqref{eq:ADMM4}, we next substitute $v^{k+1}=u^{k}+2\alpha(Ax^k-c)-(1-2\alpha)Bz^k$. Let $x_A^{k} = 2\alpha Ax^k-(1-2\alpha)(Bz^k-c)$ to get $v^{k+1}=u^{k}+x_A^k-c$ and
\begin{align*}
z^{k} &= \argmin_{z}\{g(z)+\tfrac{\rho}{2}\|x_A^{k-1}+Bz-c+u^{k-1}\|_2^2\}\\
x^k &= \argmin_{x}\{f(x)+\tfrac{\rho}{2}\|Ax+2 Bz^k+u^{k-1}+x_A^{k-1}- 2c\|_2^2\}\\
u^{k} &= u^{k-1}+(x_A^{k-1}+Bz^k-c)
\end{align*}
since  $v^{k+1}=u^{k}+x_A^k-c$ inserted in \eqref{eq:ADMM3_ls_app} implies
\begin{align*}
u^{k}&=u^{k-1}+x_A^{k-1}-x_A^k+2\alpha(Ax^{k}+Bz^k-c)\\
&=u^{k-1}+x_A^{k-1}-(2\alpha Ax^{k}-(1-2\alpha)(Bz^{k}-c))+2\alpha(Ax^{k}+Bz^k-c)\\
&=u^{k-1}+(x_A^{k-1}+Bz^k-c)
\end{align*}
(This implies that $v^{k}=u^{k}-Bz^{k}$.)
Next, insert the third equation into the second to get
\begin{align*}
z^{k} &= \argmin_{z}\{g(z)+\tfrac{\rho}{2}\|x_A^{k-1}+Bz-c+u^{k-1}\|_2^2\}\\
x^k &= \argmin_{x}\{f(x)+\tfrac{\rho}{2}\|Ax+Bz^k- c+u^{k}\|_2^2\}\\
u^{k} &= u^{k-1}+ (x_A^{k-1}+Bz^k-c)
\end{align*}
Now, change order of the $x^k$ update and the $u^{k}$ update and move the $x^{k}$ update to the first line and insert $x_A^{k-1}$ to get
\begin{align*}
x^{k-1} &= \argmin_{x}\{f(x)+\tfrac{\rho}{2}\|Ax+Bz^{k-1}-c+u^{k-1}\|_2^2\}\\
x_A^{k-1}&=2\alpha Ax^{k-1}-(1-2\alpha)(Bz^{k-1}-c)\\
z^{k} &= \argmin_{z}\{g(z)+\tfrac{\rho}{2}\|x_A^{k-1}+Bz-c+u^{k-1}\|_2^2\}\\
u^{k} &= u^{k-1}+ (x_A^{k-1}+Bz^k-c)
\end{align*}
Now, let $x^k\to x^{k+1}$ and $x_A^k\to x_A^{k+1}$ to get
\begin{align*}
x^{k} &= \argmin_{x}\{f(x)+\tfrac{\rho}{2}\|Ax+Bz^{k-1}-c+u^{k-1}\|_2^2\}\\
x_A^{k}&=2\alpha Ax^{k}-(1-2\alpha)(Bz^{k-1}-c)\\
z^{k} &= \argmin_{z}\{g(z)+\tfrac{\rho}{2}\|x_A^{k}+Bz-c+u^{k-1}\|_2^2\}\\
u^{k} &= u^{k-1}+ (x_A^{k}+Bz^{k}-c)
\end{align*}
Letting $k\to k+1$ gives ADMM on the standard form \eqref{eq:ADMM1}-\eqref{eq:ADMM4}.

\begin{rem}
ADMM can also be derived by applying Douglas-Rachford to the Fenchel dual of \eqref{eq:ADMM_prob_app}, see \cite{Gabay83}. The Fenchel dual is
\begin{align*}
{\hbox{minimize }} f^*(-A^T\mu)+c^T\mu+g^*(-B^T\mu).
\end{align*}
Letting $d_1(\mu):=f^*(-A^T\mu)+c^T\mu$ and $d_2(\mu) := g^*(-B^T\mu)$, this is equivalent to
\begin{align*}
{\hbox{minimize }} d_1(\mu)+d_2(\mu).
\end{align*}
It holds that $p_1^*=d_1$ and $p_2^*=d_2$, see \cite[Corollary~15.28]{bauschkeCVXanal}. It is also known that Douglas-Rachford when applied to minimize $p_1+p_2$ is equivalent to applying Douglas-Rachford to minimize $p_1^*+p_2^*$ (which is $d_1+d_2$), see \cite{EcksteinPhD}. Therefore we can also apply Douglas-Rachford to this dual formulation to get ADMM. This derivation is longer and therefore not used here.

\end{rem}

\end{document}

%% file: NNLS_alpha.tex
%
%
%
\definecolor{mycolor1}{rgb}{0.00000,0.44700,0.74100}%
\begin{tikzpicture}[scale=0.7]

\begin{axis}[%
width=4.6015625in,
height=3.50619791666667in,
scale only axis,
separate axis lines,
every outer x axis line/.append style={white!15!black},
every x tick label/.append style={font=\color{white!15!black}},
xmin=0,
xmax=2800,
xlabel={{\large{iteration}}},
every outer y axis line/.append style={white!15!black},
every y tick label/.append style={font=\color{white!15!black}},
ymin=0,
ymax=50,
ylabel={{\large{$\alphak$}}}
]
\addplot [color=mycolor1,solid,forget plot]
  table[row sep=crcr]{%
1	0.5\\
2	0.881928903918566\\
3	1.23470046548599\\
4	0.5\\
5	1.23470046548599\\
6	1.23470046548599\\
7	0.5\\
8	1.72858065168039\\
9	0.5\\
10	1.23470046548599\\
11	0.5\\
12	2.42001291235254\\
13	0.5\\
14	0.5\\
15	0.5\\
16	0.5\\
17	2.42001291235254\\
18	0.5\\
19	3.38801807729356\\
20	0.5\\
21	0.5\\
22	0.5\\
23	0.5\\
24	3.38801807729356\\
25	0.5\\
26	0.5\\
27	0.5\\
28	4.74322530821099\\
29	0.5\\
30	0.5\\
31	0.5\\
32	0.5\\
33	0.5\\
34	4.74322530821099\\
35	0.5\\
36	0.5\\
37	0.5\\
38	0.5\\
39	0.5\\
40	0.5\\
41	6.64051543149538\\
42	0.5\\
43	0.5\\
44	0.5\\
45	0.5\\
46	0.5\\
47	0.5\\
48	0.5\\
49	0.5\\
50	0.5\\
51	0.5\\
52	4.74322530821099\\
53	0.5\\
54	0.5\\
55	0.5\\
56	0.5\\
57	0.5\\
58	0.5\\
59	0.5\\
60	0.5\\
61	9.29672160409353\\
62	0.5\\
63	0.5\\
64	0.5\\
65	0.5\\
66	0.5\\
67	0.5\\
68	9.29672160409353\\
69	0.5\\
70	0.5\\
71	0.5\\
72	0.5\\
73	0.5\\
74	0.5\\
75	0.5\\
76	0.5\\
77	0.5\\
78	0.5\\
79	0.5\\
80	0.5\\
81	0.5\\
82	0.5\\
83	13.0154102457309\\
84	0.5\\
85	0.5\\
86	0.5\\
87	0.5\\
88	0.5\\
89	0.5\\
90	0.5\\
91	6.64051543149538\\
92	0.5\\
93	0.5\\
94	0.5\\
95	0.5\\
96	0.5\\
97	9.29672160409353\\
98	0.5\\
99	0.5\\
100	0.5\\
101	0.5\\
102	0.5\\
103	0.5\\
104	0.5\\
105	0.5\\
106	0.5\\
107	6.64051543149538\\
108	0.5\\
109	0.5\\
110	0.5\\
111	0.5\\
112	0.5\\
113	0.5\\
114	0.5\\
115	0.5\\
116	0.5\\
117	0.5\\
118	0.5\\
119	6.64051543149538\\
120	0.5\\
121	0.5\\
122	0.5\\
123	13.0154102457309\\
124	0.5\\
125	0.5\\
126	0.5\\
127	0.5\\
128	0.5\\
129	0.5\\
130	0.5\\
131	0.5\\
132	0.5\\
133	0.5\\
134	0.5\\
135	0.5\\
136	0.5\\
137	0.5\\
138	0.5\\
139	0.5\\
140	0.5\\
141	0.5\\
142	0.5\\
143	0.5\\
144	0.5\\
145	0.5\\
146	6.64051543149538\\
147	0.5\\
148	0.5\\
149	0.5\\
150	0.5\\
151	9.29672160409353\\
152	0.5\\
153	0.5\\
154	0.5\\
155	0.5\\
156	0.5\\
157	0.5\\
158	0.5\\
159	0.5\\
160	0.5\\
161	0.5\\
162	0.5\\
163	0.5\\
164	13.0154102457309\\
165	0.5\\
166	0.5\\
167	0.5\\
168	0.5\\
169	0.5\\
170	0.5\\
171	0.5\\
172	0.5\\
173	0.5\\
174	0.5\\
175	0.5\\
176	0.5\\
177	0.5\\
178	0.5\\
179	18.2215743440233\\
180	0.5\\
181	0.5\\
182	0.5\\
183	0.5\\
184	0.5\\
185	0.5\\
186	0.5\\
187	0.5\\
188	0.5\\
189	0.5\\
190	25.5102040816327\\
191	0.5\\
192	0.5\\
193	0.5\\
194	0.5\\
195	0.5\\
196	0.5\\
197	0.5\\
198	0.5\\
199	0.5\\
200	0.5\\
201	0.5\\
202	0.5\\
203	0.5\\
204	0.5\\
205	0.5\\
206	18.2215743440233\\
207	0.5\\
208	0.5\\
209	0.5\\
210	0.5\\
211	0.5\\
212	0.5\\
213	0.5\\
214	0.5\\
215	0.5\\
216	0.5\\
217	0.5\\
218	0.5\\
219	0.5\\
220	0.5\\
221	13.0154102457309\\
222	0.5\\
223	0.5\\
224	0.5\\
225	0.5\\
226	0.5\\
227	0.5\\
228	0.5\\
229	0.5\\
230	18.2215743440233\\
231	0.5\\
232	0.5\\
233	0.5\\
234	0.5\\
235	0.5\\
236	0.5\\
237	0.5\\
238	0.5\\
239	0.5\\
240	0.5\\
241	0.5\\
242	0.5\\
243	13.0154102457309\\
244	0.5\\
245	0.5\\
246	0.5\\
247	0.5\\
248	0.5\\
249	0.5\\
250	0.5\\
251	0.5\\
252	0.5\\
253	0.5\\
254	0.5\\
255	0.5\\
256	0.5\\
257	0.5\\
258	0.5\\
259	0.5\\
260	0.5\\
261	0.5\\
262	0.5\\
263	0.5\\
264	0.5\\
265	0.5\\
266	0.5\\
267	0.5\\
268	0.5\\
269	0.5\\
270	0.5\\
271	0.5\\
272	0.5\\
273	0.5\\
274	0.5\\
275	0.5\\
276	0.5\\
277	0.5\\
278	35.7142857142857\\
279	0.5\\
280	0.5\\
281	0.5\\
282	0.5\\
283	0.5\\
284	0.5\\
285	0.5\\
286	0.5\\
287	0.5\\
288	0.5\\
289	13.0154102457309\\
290	0.5\\
291	0.5\\
292	0.5\\
293	0.5\\
294	0.5\\
295	0.5\\
296	0.5\\
297	0.5\\
298	0.5\\
299	0.5\\
300	18.2215743440233\\
301	0.5\\
302	0.5\\
303	0.5\\
304	0.5\\
305	0.5\\
306	0.5\\
307	0.5\\
308	0.5\\
309	0.5\\
310	0.5\\
311	0.5\\
312	0.5\\
313	0.5\\
314	0.5\\
315	0.5\\
316	0.5\\
317	0.5\\
318	0.5\\
319	0.5\\
320	0.5\\
321	50\\
322	0.5\\
323	0.5\\
324	0.5\\
325	0.5\\
326	0.5\\
327	0.5\\
328	0.5\\
329	0.5\\
330	0.5\\
331	0.5\\
332	0.5\\
333	0.5\\
334	0.5\\
335	0.5\\
336	0.5\\
337	0.5\\
338	0.5\\
339	9.29672160409353\\
340	0.5\\
341	0.5\\
342	0.5\\
343	0.5\\
344	0.5\\
345	0.5\\
346	25.5102040816327\\
347	0.5\\
348	0.5\\
349	0.5\\
350	0.5\\
351	0.5\\
352	0.5\\
353	0.5\\
354	0.5\\
355	0.5\\
356	0.5\\
357	0.5\\
358	0.5\\
359	0.5\\
360	0.5\\
361	0.5\\
362	18.2215743440233\\
363	0.5\\
364	0.5\\
365	0.5\\
366	0.5\\
367	0.5\\
368	0.5\\
369	0.5\\
370	0.5\\
371	0.5\\
372	0.5\\
373	0.5\\
374	0.5\\
375	0.5\\
376	0.5\\
377	0.5\\
378	0.5\\
379	0.5\\
380	0.5\\
381	0.5\\
382	0.5\\
383	18.2215743440233\\
384	0.5\\
385	0.5\\
386	0.5\\
387	0.5\\
388	0.5\\
389	0.5\\
390	0.5\\
391	0.5\\
392	0.5\\
393	0.5\\
394	0.5\\
395	25.5102040816327\\
396	0.5\\
397	0.5\\
398	0.5\\
399	0.5\\
400	0.5\\
401	0.5\\
402	0.5\\
403	0.5\\
404	0.5\\
405	0.5\\
406	0.5\\
407	0.5\\
408	0.5\\
409	0.5\\
410	18.2215743440233\\
411	0.5\\
412	0.5\\
413	0.5\\
414	0.5\\
415	0.5\\
416	0.5\\
417	0.5\\
418	0.5\\
419	0.5\\
420	0.5\\
421	0.5\\
422	0.5\\
423	0.5\\
424	0.5\\
425	0.5\\
426	18.2215743440233\\
427	0.5\\
428	0.5\\
429	0.5\\
430	0.5\\
431	0.5\\
432	0.5\\
433	0.5\\
434	0.5\\
435	0.5\\
436	0.5\\
437	0.5\\
438	0.5\\
439	0.5\\
440	25.5102040816327\\
441	0.5\\
442	0.5\\
443	0.5\\
444	0.5\\
445	0.5\\
446	0.5\\
447	0.5\\
448	0.5\\
449	0.5\\
450	0.5\\
451	0.5\\
452	0.5\\
453	0.5\\
454	0.5\\
455	0.5\\
456	0.5\\
457	0.5\\
458	0.5\\
459	0.5\\
460	18.2215743440233\\
461	0.5\\
462	0.5\\
463	0.5\\
464	0.5\\
465	0.5\\
466	0.5\\
467	0.5\\
468	0.5\\
469	0.5\\
470	0.5\\
471	0.5\\
472	0.5\\
473	0.5\\
474	18.2215743440233\\
475	0.5\\
476	0.5\\
477	0.5\\
478	0.5\\
479	0.5\\
480	0.5\\
481	0.5\\
482	0.5\\
483	0.5\\
484	0.5\\
485	0.5\\
486	0.5\\
487	0.5\\
488	0.5\\
489	18.2215743440233\\
490	0.5\\
491	0.5\\
492	0.5\\
493	0.5\\
494	0.5\\
495	0.5\\
496	0.5\\
497	0.5\\
498	0.5\\
499	0.5\\
500	0.5\\
501	0.5\\
502	0.5\\
503	0.5\\
504	0.5\\
505	0.5\\
506	18.2215743440233\\
507	0.5\\
508	0.5\\
509	0.5\\
510	0.5\\
511	0.5\\
512	0.5\\
513	0.5\\
514	0.5\\
515	0.5\\
516	0.5\\
517	0.5\\
518	0.5\\
519	18.2215743440233\\
520	0.5\\
521	0.5\\
522	0.5\\
523	0.5\\
524	0.5\\
525	0.5\\
526	0.5\\
527	0.5\\
528	0.5\\
529	0.5\\
530	0.5\\
531	0.5\\
532	0.5\\
533	0.5\\
534	0.5\\
535	0.5\\
536	18.2215743440233\\
537	0.5\\
538	0.5\\
539	0.5\\
540	0.5\\
541	0.5\\
542	0.5\\
543	0.5\\
544	0.5\\
545	0.5\\
546	0.5\\
547	0.5\\
548	0.5\\
549	0.5\\
550	13.0154102457309\\
551	0.5\\
552	0.5\\
553	0.5\\
554	0.5\\
555	0.5\\
556	0.5\\
557	0.5\\
558	0.5\\
559	0.5\\
560	0.5\\
561	0.5\\
562	25.5102040816327\\
563	0.5\\
564	0.5\\
565	0.5\\
566	0.5\\
567	0.5\\
568	0.5\\
569	0.5\\
570	0.5\\
571	0.5\\
572	0.5\\
573	0.5\\
574	0.5\\
575	0.5\\
576	0.5\\
577	0.5\\
578	0.5\\
579	0.5\\
580	0.5\\
581	0.5\\
582	25.5102040816327\\
583	0.5\\
584	0.5\\
585	0.5\\
586	0.5\\
587	0.5\\
588	0.5\\
589	0.5\\
590	0.5\\
591	0.5\\
592	0.5\\
593	0.5\\
594	0.5\\
595	0.5\\
596	0.5\\
597	0.5\\
598	18.2215743440233\\
599	0.5\\
600	0.5\\
601	0.5\\
602	0.5\\
603	0.5\\
604	0.5\\
605	0.5\\
606	0.5\\
607	0.5\\
608	0.5\\
609	0.5\\
610	0.5\\
611	0.5\\
612	0.5\\
613	0.5\\
614	0.5\\
615	0.5\\
616	25.5102040816327\\
617	0.5\\
618	0.5\\
619	0.5\\
620	0.5\\
621	0.5\\
622	0.5\\
623	0.5\\
624	0.5\\
625	0.5\\
626	0.5\\
627	0.5\\
628	0.5\\
629	0.5\\
630	0.5\\
631	0.5\\
632	25.5102040816327\\
633	0.5\\
634	0.5\\
635	0.5\\
636	0.5\\
637	0.5\\
638	0.5\\
639	0.5\\
640	0.5\\
641	0.5\\
642	0.5\\
643	0.5\\
644	0.5\\
645	0.5\\
646	0.5\\
647	0.5\\
648	0.5\\
649	0.5\\
650	0.5\\
651	0.5\\
652	25.5102040816327\\
653	0.5\\
654	0.5\\
655	0.5\\
656	0.5\\
657	0.5\\
658	0.5\\
659	0.5\\
660	0.5\\
661	0.5\\
662	0.5\\
663	0.5\\
664	0.5\\
665	0.5\\
666	0.5\\
667	0.5\\
668	18.2215743440233\\
669	0.5\\
670	0.5\\
671	0.5\\
672	0.5\\
673	0.5\\
674	0.5\\
675	0.5\\
676	0.5\\
677	0.5\\
678	0.5\\
679	0.5\\
680	0.5\\
681	0.5\\
682	0.5\\
683	0.5\\
684	0.5\\
685	18.2215743440233\\
686	0.5\\
687	0.5\\
688	0.5\\
689	0.5\\
690	0.5\\
691	0.5\\
692	0.5\\
693	0.5\\
694	0.5\\
695	0.5\\
696	0.5\\
697	0.5\\
698	0.5\\
699	25.5102040816327\\
700	0.5\\
701	0.5\\
702	0.5\\
703	0.5\\
704	0.5\\
705	0.5\\
706	0.5\\
707	0.5\\
708	0.5\\
709	0.5\\
710	0.5\\
711	0.5\\
712	0.5\\
713	0.5\\
714	0.5\\
715	0.5\\
716	0.5\\
717	0.5\\
718	0.5\\
719	25.5102040816327\\
720	0.5\\
721	0.5\\
722	0.5\\
723	0.5\\
724	0.5\\
725	0.5\\
726	0.5\\
727	0.5\\
728	0.5\\
729	0.5\\
730	0.5\\
731	0.5\\
732	0.5\\
733	0.5\\
734	0.5\\
735	25.5102040816327\\
736	0.5\\
737	0.5\\
738	0.5\\
739	0.5\\
740	0.5\\
741	0.5\\
742	0.5\\
743	0.5\\
744	0.5\\
745	0.5\\
746	0.5\\
747	0.5\\
748	0.5\\
749	0.5\\
750	0.5\\
751	0.5\\
752	0.5\\
753	0.5\\
754	0.5\\
755	25.5102040816327\\
756	0.5\\
757	0.5\\
758	0.5\\
759	0.5\\
760	0.5\\
761	0.5\\
762	0.5\\
763	0.5\\
764	0.5\\
765	0.5\\
766	0.5\\
767	0.5\\
768	0.5\\
769	0.5\\
770	0.5\\
771	18.2215743440233\\
772	0.5\\
773	0.5\\
774	0.5\\
775	0.5\\
776	0.5\\
777	0.5\\
778	0.5\\
779	0.5\\
780	0.5\\
781	0.5\\
782	0.5\\
783	0.5\\
784	0.5\\
785	0.5\\
786	0.5\\
787	0.5\\
788	18.2215743440233\\
789	0.5\\
790	0.5\\
791	0.5\\
792	0.5\\
793	0.5\\
794	0.5\\
795	0.5\\
796	0.5\\
797	0.5\\
798	0.5\\
799	0.5\\
800	0.5\\
801	0.5\\
802	18.2215743440233\\
803	0.5\\
804	0.5\\
805	0.5\\
806	0.5\\
807	0.5\\
808	0.5\\
809	0.5\\
810	0.5\\
811	0.5\\
812	0.5\\
813	0.5\\
814	0.5\\
815	0.5\\
816	0.5\\
817	0.5\\
818	0.5\\
819	13.0154102457309\\
820	0.5\\
821	0.5\\
822	0.5\\
823	0.5\\
824	0.5\\
825	0.5\\
826	0.5\\
827	0.5\\
828	0.5\\
829	0.5\\
830	0.5\\
831	25.5102040816327\\
832	0.5\\
833	0.5\\
834	0.5\\
835	0.5\\
836	0.5\\
837	0.5\\
838	0.5\\
839	0.5\\
840	0.5\\
841	0.5\\
842	0.5\\
843	0.5\\
844	0.5\\
845	0.5\\
846	0.5\\
847	18.2215743440233\\
848	0.5\\
849	0.5\\
850	0.5\\
851	0.5\\
852	0.5\\
853	0.5\\
854	0.5\\
855	0.5\\
856	0.5\\
857	0.5\\
858	0.5\\
859	0.5\\
860	0.5\\
861	0.5\\
862	0.5\\
863	0.5\\
864	18.2215743440233\\
865	0.5\\
866	0.5\\
867	0.5\\
868	0.5\\
869	0.5\\
870	0.5\\
871	0.5\\
872	0.5\\
873	0.5\\
874	0.5\\
875	0.5\\
876	0.5\\
877	0.5\\
878	0.5\\
879	25.5102040816327\\
880	0.5\\
881	0.5\\
882	0.5\\
883	0.5\\
884	0.5\\
885	0.5\\
886	0.5\\
887	0.5\\
888	0.5\\
889	0.5\\
890	0.5\\
891	0.5\\
892	0.5\\
893	0.5\\
894	0.5\\
895	0.5\\
896	0.5\\
897	18.2215743440233\\
898	0.5\\
899	0.5\\
900	0.5\\
901	0.5\\
902	0.5\\
903	0.5\\
904	0.5\\
905	0.5\\
906	0.5\\
907	0.5\\
908	0.5\\
909	0.5\\
910	0.5\\
911	0.5\\
912	0.5\\
913	35.7142857142857\\
914	2.42001291235254\\
915	0.5\\
916	0.5\\
917	0.5\\
918	0.5\\
919	0.5\\
920	0.5\\
921	0.5\\
922	0.5\\
923	0.5\\
924	0.5\\
925	0.5\\
926	0.5\\
927	0.5\\
928	0.5\\
929	0.5\\
930	0.5\\
931	0.5\\
932	0.5\\
933	35.7142857142857\\
934	2.42001291235254\\
935	0.5\\
936	0.5\\
937	0.5\\
938	0.5\\
939	0.5\\
940	0.5\\
941	0.5\\
942	0.5\\
943	0.5\\
944	0.5\\
945	0.5\\
946	0.5\\
947	0.5\\
948	0.5\\
949	0.5\\
950	0.5\\
951	0.5\\
952	18.2215743440233\\
953	0.5\\
954	0.5\\
955	0.5\\
956	0.5\\
957	0.5\\
958	0.5\\
959	0.5\\
960	0.5\\
961	0.5\\
962	0.5\\
963	0.5\\
964	0.5\\
965	0.5\\
966	25.5102040816327\\
967	0.5\\
968	0.5\\
969	0.5\\
970	0.5\\
971	0.5\\
972	0.5\\
973	0.5\\
974	0.5\\
975	0.5\\
976	0.5\\
977	0.5\\
978	0.5\\
979	0.5\\
980	0.5\\
981	0.5\\
982	0.5\\
983	0.5\\
984	0.5\\
985	18.2215743440233\\
986	0.5\\
987	0.5\\
988	0.5\\
989	0.5\\
990	0.5\\
991	0.5\\
992	0.5\\
993	0.5\\
994	0.5\\
995	0.5\\
996	0.5\\
997	0.5\\
998	0.5\\
999	0.5\\
1000	25.5102040816327\\
1001	0.5\\
1002	0.5\\
1003	0.5\\
1004	0.5\\
1005	0.5\\
1006	0.5\\
1007	0.5\\
1008	0.5\\
1009	0.5\\
1010	0.5\\
1011	0.5\\
1012	0.5\\
1013	0.5\\
1014	0.5\\
1015	0.5\\
1016	0.5\\
1017	0.5\\
1018	25.5102040816327\\
1019	0.5\\
1020	0.5\\
1021	0.5\\
1022	0.5\\
1023	0.5\\
1024	0.5\\
1025	0.5\\
1026	0.5\\
1027	0.5\\
1028	0.5\\
1029	0.5\\
1030	0.5\\
1031	0.5\\
1032	0.5\\
1033	0.5\\
1034	0.5\\
1035	0.5\\
1036	25.5102040816327\\
1037	0.5\\
1038	0.5\\
1039	0.5\\
1040	0.5\\
1041	0.5\\
1042	0.5\\
1043	0.5\\
1044	0.5\\
1045	0.5\\
1046	0.5\\
1047	0.5\\
1048	0.5\\
1049	0.5\\
1050	0.5\\
1051	0.5\\
1052	0.5\\
1053	0.5\\
1054	0.5\\
1055	25.5102040816327\\
1056	0.5\\
1057	0.5\\
1058	0.5\\
1059	0.5\\
1060	0.5\\
1061	0.5\\
1062	0.5\\
1063	0.5\\
1064	0.5\\
1065	0.5\\
1066	0.5\\
1067	0.5\\
1068	0.5\\
1069	0.5\\
1070	0.5\\
1071	0.5\\
1072	18.2215743440233\\
1073	0.5\\
1074	0.5\\
1075	0.5\\
1076	0.5\\
1077	0.5\\
1078	0.5\\
1079	0.5\\
1080	0.5\\
1081	0.5\\
1082	0.5\\
1083	0.5\\
1084	0.5\\
1085	0.5\\
1086	0.5\\
1087	18.2215743440233\\
1088	0.5\\
1089	0.5\\
1090	0.5\\
1091	0.5\\
1092	0.5\\
1093	0.5\\
1094	0.5\\
1095	0.5\\
1096	0.5\\
1097	0.5\\
1098	0.5\\
1099	0.5\\
1100	0.5\\
1101	0.5\\
1102	0.5\\
1103	0.5\\
1104	25.5102040816327\\
1105	0.5\\
1106	0.5\\
1107	0.5\\
1108	0.5\\
1109	0.5\\
1110	0.5\\
1111	0.5\\
1112	0.5\\
1113	0.5\\
1114	0.5\\
1115	0.5\\
1116	0.5\\
1117	0.5\\
1118	0.5\\
1119	0.5\\
1120	25.5102040816327\\
1121	0.5\\
1122	0.5\\
1123	0.5\\
1124	0.5\\
1125	0.5\\
1126	0.5\\
1127	0.5\\
1128	0.5\\
1129	0.5\\
1130	0.5\\
1131	0.5\\
1132	0.5\\
1133	0.5\\
1134	0.5\\
1135	0.5\\
1136	0.5\\
1137	0.5\\
1138	0.5\\
1139	0.5\\
1140	25.5102040816327\\
1141	0.5\\
1142	0.5\\
1143	0.5\\
1144	0.5\\
1145	0.5\\
1146	0.5\\
1147	0.5\\
1148	0.5\\
1149	0.5\\
1150	0.5\\
1151	0.5\\
1152	0.5\\
1153	0.5\\
1154	0.5\\
1155	0.5\\
1156	0.5\\
1157	25.5102040816327\\
1158	0.5\\
1159	0.5\\
1160	0.5\\
1161	0.5\\
1162	0.5\\
1163	0.5\\
1164	0.5\\
1165	0.5\\
1166	0.5\\
1167	0.5\\
1168	0.5\\
1169	0.5\\
1170	0.5\\
1171	0.5\\
1172	0.5\\
1173	0.5\\
1174	0.5\\
1175	0.5\\
1176	25.5102040816327\\
1177	0.5\\
1178	0.5\\
1179	0.5\\
1180	0.5\\
1181	0.5\\
1182	0.5\\
1183	0.5\\
1184	0.5\\
1185	0.5\\
1186	0.5\\
1187	0.5\\
1188	0.5\\
1189	0.5\\
1190	0.5\\
1191	0.5\\
1192	0.5\\
1193	25.5102040816327\\
1194	0.5\\
1195	0.5\\
1196	0.5\\
1197	0.5\\
1198	0.5\\
1199	0.5\\
1200	0.5\\
1201	0.5\\
1202	0.5\\
1203	0.5\\
1204	0.5\\
1205	0.5\\
1206	0.5\\
1207	0.5\\
1208	0.5\\
1209	0.5\\
1210	0.5\\
1211	0.5\\
1212	25.5102040816327\\
1213	0.5\\
1214	0.5\\
1215	0.5\\
1216	0.5\\
1217	0.5\\
1218	0.5\\
1219	0.5\\
1220	0.5\\
1221	0.5\\
1222	0.5\\
1223	0.5\\
1224	0.5\\
1225	0.5\\
1226	0.5\\
1227	0.5\\
1228	0.5\\
1229	0.5\\
1230	25.5102040816327\\
1231	0.5\\
1232	0.5\\
1233	0.5\\
1234	0.5\\
1235	0.5\\
1236	0.5\\
1237	0.5\\
1238	0.5\\
1239	0.5\\
1240	0.5\\
1241	0.5\\
1242	0.5\\
1243	0.5\\
1244	0.5\\
1245	0.5\\
1246	0.5\\
1247	0.5\\
1248	25.5102040816327\\
1249	0.5\\
1250	0.5\\
1251	0.5\\
1252	0.5\\
1253	0.5\\
1254	0.5\\
1255	0.5\\
1256	0.5\\
1257	0.5\\
1258	0.5\\
1259	0.5\\
1260	0.5\\
1261	0.5\\
1262	0.5\\
1263	0.5\\
1264	0.5\\
1265	0.5\\
1266	18.2215743440233\\
1267	0.5\\
1268	0.5\\
1269	0.5\\
1270	0.5\\
1271	0.5\\
1272	0.5\\
1273	0.5\\
1274	0.5\\
1275	0.5\\
1276	0.5\\
1277	0.5\\
1278	0.5\\
1279	0.5\\
1280	0.5\\
1281	0.5\\
1282	25.5102040816327\\
1283	0.5\\
1284	0.5\\
1285	0.5\\
1286	0.5\\
1287	0.5\\
1288	0.5\\
1289	0.5\\
1290	0.5\\
1291	0.5\\
1292	0.5\\
1293	0.5\\
1294	0.5\\
1295	0.5\\
1296	0.5\\
1297	0.5\\
1298	0.5\\
1299	25.5102040816327\\
1300	0.5\\
1301	0.5\\
1302	0.5\\
1303	0.5\\
1304	0.5\\
1305	0.5\\
1306	0.5\\
1307	0.5\\
1308	0.5\\
1309	0.5\\
1310	0.5\\
1311	0.5\\
1312	0.5\\
1313	0.5\\
1314	0.5\\
1315	0.5\\
1316	0.5\\
1317	0.5\\
1318	25.5102040816327\\
1319	0.5\\
1320	0.5\\
1321	0.5\\
1322	0.5\\
1323	0.5\\
1324	0.5\\
1325	0.5\\
1326	0.5\\
1327	0.5\\
1328	0.5\\
1329	0.5\\
1330	0.5\\
1331	0.5\\
1332	0.5\\
1333	0.5\\
1334	0.5\\
1335	18.2215743440233\\
1336	0.5\\
1337	0.5\\
1338	0.5\\
1339	0.5\\
1340	0.5\\
1341	0.5\\
1342	0.5\\
1343	0.5\\
1344	0.5\\
1345	0.5\\
1346	0.5\\
1347	0.5\\
1348	0.5\\
1349	0.5\\
1350	0.5\\
1351	0.5\\
1352	35.7142857142857\\
1353	2.42001291235254\\
1354	0.5\\
1355	0.5\\
1356	0.5\\
1357	0.5\\
1358	0.5\\
1359	0.5\\
1360	0.5\\
1361	0.5\\
1362	0.5\\
1363	0.5\\
1364	0.5\\
1365	0.5\\
1366	0.5\\
1367	0.5\\
1368	0.5\\
1369	0.5\\
1370	0.5\\
1371	0.5\\
1372	35.7142857142857\\
1373	2.42001291235254\\
1374	0.5\\
1375	0.5\\
1376	0.5\\
1377	0.5\\
1378	0.5\\
1379	0.5\\
1380	0.5\\
1381	0.5\\
1382	0.5\\
1383	0.5\\
1384	0.5\\
1385	0.5\\
1386	0.5\\
1387	0.5\\
1388	0.5\\
1389	0.5\\
1390	0.5\\
1391	25.5102040816327\\
1392	0.5\\
1393	0.5\\
1394	0.5\\
1395	0.5\\
1396	0.5\\
1397	0.5\\
1398	0.5\\
1399	0.5\\
1400	0.5\\
1401	0.5\\
1402	0.5\\
1403	0.5\\
1404	0.5\\
1405	0.5\\
1406	0.5\\
1407	25.5102040816327\\
1408	0.5\\
1409	0.5\\
1410	0.5\\
1411	0.5\\
1412	0.5\\
1413	0.5\\
1414	0.5\\
1415	0.5\\
1416	0.5\\
1417	0.5\\
1418	0.5\\
1419	0.5\\
1420	0.5\\
1421	0.5\\
1422	0.5\\
1423	0.5\\
1424	0.5\\
1425	0.5\\
1426	0.5\\
1427	25.5102040816327\\
1428	0.5\\
1429	0.5\\
1430	0.5\\
1431	0.5\\
1432	0.5\\
1433	0.5\\
1434	0.5\\
1435	0.5\\
1436	0.5\\
1437	0.5\\
1438	0.5\\
1439	0.5\\
1440	0.5\\
1441	0.5\\
1442	0.5\\
1443	0.5\\
1444	35.7142857142857\\
1445	0.5\\
1446	0.5\\
1447	0.5\\
1448	0.5\\
1449	0.5\\
1450	0.5\\
1451	0.5\\
1452	0.5\\
1453	0.5\\
1454	0.5\\
1455	0.5\\
1456	0.5\\
1457	0.5\\
1458	0.5\\
1459	0.5\\
1460	0.5\\
1461	0.5\\
1462	0.5\\
1463	0.5\\
1464	0.5\\
1465	0.5\\
1466	25.5102040816327\\
1467	0.5\\
1468	0.5\\
1469	0.5\\
1470	0.5\\
1471	0.5\\
1472	0.5\\
1473	0.5\\
1474	0.5\\
1475	0.5\\
1476	0.5\\
1477	0.5\\
1478	0.5\\
1479	0.5\\
1480	0.5\\
1481	0.5\\
1482	25.5102040816327\\
1483	0.5\\
1484	0.5\\
1485	0.5\\
1486	0.5\\
1487	0.5\\
1488	0.5\\
1489	0.5\\
1490	0.5\\
1491	0.5\\
1492	0.5\\
1493	0.5\\
1494	0.5\\
1495	0.5\\
1496	0.5\\
1497	0.5\\
1498	0.5\\
1499	0.5\\
1500	0.5\\
1501	0.5\\
1502	25.5102040816327\\
1503	0.5\\
1504	0.5\\
1505	0.5\\
1506	0.5\\
1507	0.5\\
1508	0.5\\
1509	0.5\\
1510	0.5\\
1511	0.5\\
1512	0.5\\
1513	0.5\\
1514	0.5\\
1515	0.5\\
1516	0.5\\
1517	0.5\\
1518	0.5\\
1519	25.5102040816327\\
1520	0.5\\
1521	0.5\\
1522	0.5\\
1523	0.5\\
1524	0.5\\
1525	0.5\\
1526	0.5\\
1527	0.5\\
1528	0.5\\
1529	0.5\\
1530	0.5\\
1531	0.5\\
1532	0.5\\
1533	0.5\\
1534	0.5\\
1535	0.5\\
1536	0.5\\
1537	0.5\\
1538	25.5102040816327\\
1539	0.5\\
1540	0.5\\
1541	0.5\\
1542	0.5\\
1543	0.5\\
1544	0.5\\
1545	0.5\\
1546	0.5\\
1547	0.5\\
1548	0.5\\
1549	0.5\\
1550	0.5\\
1551	0.5\\
1552	0.5\\
1553	0.5\\
1554	0.5\\
1555	25.5102040816327\\
1556	0.5\\
1557	0.5\\
1558	0.5\\
1559	0.5\\
1560	0.5\\
1561	0.5\\
1562	0.5\\
1563	0.5\\
1564	0.5\\
1565	0.5\\
1566	0.5\\
1567	0.5\\
1568	0.5\\
1569	0.5\\
1570	0.5\\
1571	0.5\\
1572	0.5\\
1573	0.5\\
1574	25.5102040816327\\
1575	0.5\\
1576	0.5\\
1577	0.5\\
1578	0.5\\
1579	0.5\\
1580	0.5\\
1581	0.5\\
1582	0.5\\
1583	0.5\\
1584	0.5\\
1585	0.5\\
1586	0.5\\
1587	0.5\\
1588	0.5\\
1589	0.5\\
1590	0.5\\
1591	0.5\\
1592	25.5102040816327\\
1593	0.5\\
1594	0.5\\
1595	0.5\\
1596	0.5\\
1597	0.5\\
1598	0.5\\
1599	0.5\\
1600	0.5\\
1601	0.5\\
1602	0.5\\
1603	0.5\\
1604	0.5\\
1605	0.5\\
1606	0.5\\
1607	0.5\\
1608	0.5\\
1609	0.5\\
1610	25.5102040816327\\
1611	0.5\\
1612	0.5\\
1613	0.5\\
1614	0.5\\
1615	0.5\\
1616	0.5\\
1617	0.5\\
1618	0.5\\
1619	0.5\\
1620	0.5\\
1621	0.5\\
1622	0.5\\
1623	0.5\\
1624	0.5\\
1625	0.5\\
1626	0.5\\
1627	0.5\\
1628	0.5\\
1629	35.7142857142857\\
1630	0.5\\
1631	0.5\\
1632	0.5\\
1633	0.5\\
1634	0.5\\
1635	0.5\\
1636	0.5\\
1637	0.5\\
1638	0.5\\
1639	0.5\\
1640	0.5\\
1641	0.5\\
1642	0.5\\
1643	0.5\\
1644	0.5\\
1645	0.5\\
1646	0.5\\
1647	0.5\\
1648	25.5102040816327\\
1649	0.5\\
1650	0.5\\
1651	0.5\\
1652	0.5\\
1653	0.5\\
1654	0.5\\
1655	0.5\\
1656	0.5\\
1657	0.5\\
1658	0.5\\
1659	0.5\\
1660	0.5\\
1661	0.5\\
1662	0.5\\
1663	0.5\\
1664	0.5\\
1665	0.5\\
1666	0.5\\
1667	0.5\\
1668	25.5102040816327\\
1669	0.5\\
1670	0.5\\
1671	0.5\\
1672	0.5\\
1673	0.5\\
1674	0.5\\
1675	0.5\\
1676	0.5\\
1677	0.5\\
1678	0.5\\
1679	0.5\\
1680	0.5\\
1681	0.5\\
1682	0.5\\
1683	0.5\\
1684	25.5102040816327\\
1685	0.5\\
1686	0.5\\
1687	0.5\\
1688	0.5\\
1689	0.5\\
1690	0.5\\
1691	0.5\\
1692	0.5\\
1693	0.5\\
1694	0.5\\
1695	0.5\\
1696	0.5\\
1697	0.5\\
1698	0.5\\
1699	0.5\\
1700	0.5\\
1701	0.5\\
1702	0.5\\
1703	0.5\\
1704	25.5102040816327\\
1705	0.5\\
1706	0.5\\
1707	0.5\\
1708	0.5\\
1709	0.5\\
1710	0.5\\
1711	0.5\\
1712	0.5\\
1713	0.5\\
1714	0.5\\
1715	0.5\\
1716	0.5\\
1717	0.5\\
1718	0.5\\
1719	0.5\\
1720	25.5102040816327\\
1721	0.5\\
1722	0.5\\
1723	0.5\\
1724	0.5\\
1725	0.5\\
1726	0.5\\
1727	0.5\\
1728	0.5\\
1729	0.5\\
1730	0.5\\
1731	0.5\\
1732	0.5\\
1733	0.5\\
1734	0.5\\
1735	0.5\\
1736	0.5\\
1737	0.5\\
1738	0.5\\
1739	0.5\\
1740	25.5102040816327\\
1741	0.5\\
1742	0.5\\
1743	0.5\\
1744	0.5\\
1745	0.5\\
1746	0.5\\
1747	0.5\\
1748	0.5\\
1749	0.5\\
1750	0.5\\
1751	0.5\\
1752	0.5\\
1753	0.5\\
1754	0.5\\
1755	0.5\\
1756	25.5102040816327\\
1757	0.5\\
1758	0.5\\
1759	0.5\\
1760	0.5\\
1761	0.5\\
1762	0.5\\
1763	0.5\\
1764	0.5\\
1765	0.5\\
1766	0.5\\
1767	0.5\\
1768	0.5\\
1769	0.5\\
1770	0.5\\
1771	0.5\\
1772	0.5\\
1773	0.5\\
1774	0.5\\
1775	0.5\\
1776	0.5\\
1777	35.7142857142857\\
1778	0.5\\
1779	0.5\\
1780	0.5\\
1781	0.5\\
1782	0.5\\
1783	0.5\\
1784	0.5\\
1785	0.5\\
1786	0.5\\
1787	0.5\\
1788	0.5\\
1789	0.5\\
1790	0.5\\
1791	0.5\\
1792	0.5\\
1793	0.5\\
1794	0.5\\
1795	0.5\\
1796	0.5\\
1797	25.5102040816327\\
1798	0.5\\
1799	0.5\\
1800	0.5\\
1801	0.5\\
1802	0.5\\
1803	0.5\\
1804	0.5\\
1805	0.5\\
1806	0.5\\
1807	0.5\\
1808	0.5\\
1809	0.5\\
1810	0.5\\
1811	0.5\\
1812	0.5\\
1813	0.5\\
1814	0.5\\
1815	35.7142857142857\\
1816	0.5\\
1817	0.5\\
1818	0.5\\
1819	0.5\\
1820	0.5\\
1821	0.5\\
1822	0.5\\
1823	0.5\\
1824	0.5\\
1825	0.5\\
1826	0.5\\
1827	0.5\\
1828	0.5\\
1829	0.5\\
1830	0.5\\
1831	0.5\\
1832	0.5\\
1833	0.5\\
1834	0.5\\
1835	0.5\\
1836	0.5\\
1837	35.7142857142857\\
1838	0.5\\
1839	0.5\\
1840	0.5\\
1841	0.5\\
1842	0.5\\
1843	0.5\\
1844	0.5\\
1845	0.5\\
1846	0.5\\
1847	0.5\\
1848	0.5\\
1849	0.5\\
1850	0.5\\
1851	0.5\\
1852	0.5\\
1853	0.5\\
1854	0.5\\
1855	0.5\\
1856	35.7142857142857\\
1857	0.5\\
1858	0.5\\
1859	0.5\\
1860	0.5\\
1861	0.5\\
1862	0.5\\
1863	0.5\\
1864	0.5\\
1865	0.5\\
1866	0.5\\
1867	0.5\\
1868	0.5\\
1869	0.5\\
1870	0.5\\
1871	0.5\\
1872	0.5\\
1873	0.5\\
1874	0.5\\
1875	0.5\\
1876	0.5\\
1877	0.5\\
1878	35.7142857142857\\
1879	0.5\\
1880	0.5\\
1881	0.5\\
1882	0.5\\
1883	0.5\\
1884	0.5\\
1885	0.5\\
1886	0.5\\
1887	0.5\\
1888	0.5\\
1889	0.5\\
1890	0.5\\
1891	0.5\\
1892	0.5\\
1893	0.5\\
1894	0.5\\
1895	0.5\\
1896	0.5\\
1897	0.5\\
1898	0.5\\
1899	0.5\\
1900	35.7142857142857\\
1901	0.5\\
1902	0.5\\
1903	0.5\\
1904	0.5\\
1905	0.5\\
1906	0.5\\
1907	0.5\\
1908	0.5\\
1909	0.5\\
1910	0.5\\
1911	0.5\\
1912	0.5\\
1913	0.5\\
1914	0.5\\
1915	0.5\\
1916	0.5\\
1917	0.5\\
1918	0.5\\
1919	35.7142857142857\\
1920	0.5\\
1921	0.5\\
1922	0.5\\
1923	0.5\\
1924	0.5\\
1925	0.5\\
1926	0.5\\
1927	0.5\\
1928	0.5\\
1929	0.5\\
1930	0.5\\
1931	0.5\\
1932	0.5\\
1933	0.5\\
1934	0.5\\
1935	0.5\\
1936	0.5\\
1937	0.5\\
1938	0.5\\
1939	0.5\\
1940	0.5\\
1941	35.7142857142857\\
1942	0.5\\
1943	0.5\\
1944	0.5\\
1945	0.5\\
1946	0.5\\
1947	0.5\\
1948	0.5\\
1949	0.5\\
1950	0.5\\
1951	0.5\\
1952	0.5\\
1953	0.5\\
1954	0.5\\
1955	0.5\\
1956	0.5\\
1957	0.5\\
1958	0.5\\
1959	0.5\\
1960	0.5\\
1961	0.5\\
1962	0.5\\
1963	35.7142857142857\\
1964	0.5\\
1965	0.5\\
1966	0.5\\
1967	0.5\\
1968	0.5\\
1969	0.5\\
1970	0.5\\
1971	0.5\\
1972	0.5\\
1973	0.5\\
1974	0.5\\
1975	0.5\\
1976	0.5\\
1977	0.5\\
1978	0.5\\
1979	0.5\\
1980	0.5\\
1981	0.5\\
1982	35.7142857142857\\
1983	0.5\\
1984	0.5\\
1985	0.5\\
1986	0.5\\
1987	0.5\\
1988	0.5\\
1989	0.5\\
1990	0.5\\
1991	0.5\\
1992	0.5\\
1993	0.5\\
1994	0.5\\
1995	0.5\\
1996	0.5\\
1997	0.5\\
1998	0.5\\
1999	0.5\\
2000	0.5\\
2001	0.5\\
2002	0.5\\
2003	0.5\\
2004	35.7142857142857\\
2005	0.5\\
2006	0.5\\
2007	0.5\\
2008	0.5\\
2009	0.5\\
2010	0.5\\
2011	0.5\\
2012	0.5\\
2013	0.5\\
2014	0.5\\
2015	0.5\\
2016	0.5\\
2017	0.5\\
2018	0.5\\
2019	0.5\\
2020	0.5\\
2021	0.5\\
2022	0.5\\
2023	0.5\\
2024	0.5\\
2025	0.5\\
2026	35.7142857142857\\
2027	0.5\\
2028	0.5\\
2029	0.5\\
2030	0.5\\
2031	0.5\\
2032	0.5\\
2033	0.5\\
2034	0.5\\
2035	0.5\\
2036	0.5\\
2037	0.5\\
2038	0.5\\
2039	0.5\\
2040	0.5\\
2041	0.5\\
2042	0.5\\
2043	0.5\\
2044	0.5\\
2045	35.7142857142857\\
2046	0.5\\
2047	0.5\\
2048	0.5\\
2049	0.5\\
2050	0.5\\
2051	0.5\\
2052	0.5\\
2053	0.5\\
2054	0.5\\
2055	0.5\\
2056	0.5\\
2057	0.5\\
2058	0.5\\
2059	0.5\\
2060	0.5\\
2061	0.5\\
2062	0.5\\
2063	0.5\\
2064	0.5\\
2065	0.5\\
2066	0.5\\
2067	35.7142857142857\\
2068	0.5\\
2069	0.5\\
2070	0.5\\
2071	0.5\\
2072	0.5\\
2073	0.5\\
2074	0.5\\
2075	0.5\\
2076	0.5\\
2077	0.5\\
2078	0.5\\
2079	0.5\\
2080	0.5\\
2081	0.5\\
2082	0.5\\
2083	0.5\\
2084	0.5\\
2085	0.5\\
2086	0.5\\
2087	0.5\\
2088	25.5102040816327\\
2089	0.5\\
2090	0.5\\
2091	0.5\\
2092	0.5\\
2093	0.5\\
2094	0.5\\
2095	0.5\\
2096	0.5\\
2097	0.5\\
2098	0.5\\
2099	0.5\\
2100	0.5\\
2101	0.5\\
2102	0.5\\
2103	0.5\\
2104	0.5\\
2105	25.5102040816327\\
2106	0.5\\
2107	0.5\\
2108	0.5\\
2109	0.5\\
2110	0.5\\
2111	0.5\\
2112	0.5\\
2113	0.5\\
2114	0.5\\
2115	0.5\\
2116	0.5\\
2117	0.5\\
2118	0.5\\
2119	0.5\\
2120	0.5\\
2121	0.5\\
2122	0.5\\
2123	0.5\\
2124	0.5\\
2125	35.7142857142857\\
2126	0.5\\
2127	0.5\\
2128	0.5\\
2129	0.5\\
2130	0.5\\
2131	0.5\\
2132	0.5\\
2133	0.5\\
2134	0.5\\
2135	0.5\\
2136	0.5\\
2137	0.5\\
2138	0.5\\
2139	0.5\\
2140	0.5\\
2141	0.5\\
2142	0.5\\
2143	0.5\\
2144	35.7142857142857\\
2145	0.5\\
2146	0.5\\
2147	0.5\\
2148	0.5\\
2149	0.5\\
2150	0.5\\
2151	0.5\\
2152	0.5\\
2153	0.5\\
2154	0.5\\
2155	0.5\\
2156	0.5\\
2157	0.5\\
2158	0.5\\
2159	0.5\\
2160	0.5\\
2161	0.5\\
2162	0.5\\
2163	0.5\\
2164	0.5\\
2165	0.5\\
2166	35.7142857142857\\
2167	0.5\\
2168	0.5\\
2169	0.5\\
2170	0.5\\
2171	0.5\\
2172	0.5\\
2173	0.5\\
2174	0.5\\
2175	0.5\\
2176	0.5\\
2177	0.5\\
2178	0.5\\
2179	0.5\\
2180	0.5\\
2181	0.5\\
2182	0.5\\
2183	0.5\\
2184	0.5\\
2185	0.5\\
2186	0.5\\
2187	0.5\\
2188	35.7142857142857\\
2189	0.5\\
2190	0.5\\
2191	0.5\\
2192	0.5\\
2193	0.5\\
2194	0.5\\
2195	0.5\\
2196	0.5\\
2197	0.5\\
2198	0.5\\
2199	0.5\\
2200	0.5\\
2201	0.5\\
2202	0.5\\
2203	0.5\\
2204	0.5\\
2205	0.5\\
2206	0.5\\
2207	35.7142857142857\\
2208	0.5\\
2209	0.5\\
2210	0.5\\
2211	0.5\\
2212	0.5\\
2213	0.5\\
2214	0.5\\
2215	0.5\\
2216	0.5\\
2217	0.5\\
2218	0.5\\
2219	0.5\\
2220	0.5\\
2221	0.5\\
2222	0.5\\
2223	0.5\\
2224	0.5\\
2225	0.5\\
2226	0.5\\
2227	0.5\\
2228	0.5\\
2229	35.7142857142857\\
2230	0.5\\
2231	0.5\\
2232	0.5\\
2233	0.5\\
2234	0.5\\
2235	0.5\\
2236	0.5\\
2237	0.5\\
2238	0.5\\
2239	0.5\\
2240	0.5\\
2241	0.5\\
2242	0.5\\
2243	0.5\\
2244	0.5\\
2245	0.5\\
2246	0.5\\
2247	0.5\\
2248	0.5\\
2249	0.5\\
2250	0.5\\
2251	35.7142857142857\\
2252	0.5\\
2253	0.5\\
2254	0.5\\
2255	0.5\\
2256	0.5\\
2257	0.5\\
2258	0.5\\
2259	0.5\\
2260	0.5\\
2261	0.5\\
2262	0.5\\
2263	0.5\\
2264	0.5\\
2265	0.5\\
2266	0.5\\
2267	0.5\\
2268	0.5\\
2269	0.5\\
2270	35.7142857142857\\
2271	0.5\\
2272	0.5\\
2273	0.5\\
2274	0.5\\
2275	0.5\\
2276	0.5\\
2277	0.5\\
2278	0.5\\
2279	0.5\\
2280	0.5\\
2281	0.5\\
2282	0.5\\
2283	0.5\\
2284	0.5\\
2285	0.5\\
2286	0.5\\
2287	0.5\\
2288	0.5\\
2289	0.5\\
2290	0.5\\
2291	0.5\\
2292	35.7142857142857\\
2293	0.5\\
2294	0.5\\
2295	0.5\\
2296	0.5\\
2297	0.5\\
2298	0.5\\
2299	0.5\\
2300	0.5\\
2301	0.5\\
2302	0.5\\
2303	0.5\\
2304	0.5\\
2305	0.5\\
2306	0.5\\
2307	0.5\\
2308	0.5\\
2309	0.5\\
2310	0.5\\
2311	0.5\\
2312	0.5\\
2313	25.5102040816327\\
2314	0.5\\
2315	0.5\\
2316	0.5\\
2317	0.5\\
2318	0.5\\
2319	0.5\\
2320	0.5\\
2321	0.5\\
2322	0.5\\
2323	0.5\\
2324	0.5\\
2325	0.5\\
2326	0.5\\
2327	0.5\\
2328	0.5\\
2329	0.5\\
2330	35.7142857142857\\
2331	0.5\\
2332	0.5\\
2333	0.5\\
2334	0.5\\
2335	0.5\\
2336	0.5\\
2337	0.5\\
2338	0.5\\
2339	0.5\\
2340	0.5\\
2341	0.5\\
2342	0.5\\
2343	0.5\\
2344	0.5\\
2345	0.5\\
2346	0.5\\
2347	0.5\\
2348	0.5\\
2349	0.5\\
2350	0.5\\
2351	0.5\\
2352	0.5\\
2353	35.7142857142857\\
2354	0.5\\
2355	0.5\\
2356	0.5\\
2357	0.5\\
2358	0.5\\
2359	0.5\\
2360	0.5\\
2361	0.5\\
2362	0.5\\
2363	0.5\\
2364	0.5\\
2365	0.5\\
2366	0.5\\
2367	0.5\\
2368	0.5\\
2369	0.5\\
2370	0.5\\
2371	0.5\\
2372	0.5\\
2373	35.7142857142857\\
2374	0.5\\
2375	0.5\\
2376	0.5\\
2377	0.5\\
2378	0.5\\
2379	0.5\\
2380	0.5\\
2381	0.5\\
2382	0.5\\
2383	0.5\\
2384	0.5\\
2385	0.5\\
2386	0.5\\
2387	0.5\\
2388	0.5\\
2389	0.5\\
2390	0.5\\
2391	0.5\\
2392	0.5\\
2393	35.7142857142857\\
2394	0.5\\
2395	0.5\\
2396	0.5\\
2397	0.5\\
2398	0.5\\
2399	0.5\\
2400	0.5\\
2401	0.5\\
2402	0.5\\
2403	0.5\\
2404	0.5\\
2405	0.5\\
2406	0.5\\
2407	0.5\\
2408	0.5\\
2409	0.5\\
2410	0.5\\
2411	0.5\\
2412	0.5\\
2413	0.5\\
2414	0.5\\
2415	25.5102040816327\\
2416	0.5\\
2417	0.5\\
2418	0.5\\
2419	0.5\\
2420	0.5\\
2421	0.5\\
2422	0.5\\
2423	0.5\\
2424	0.5\\
2425	0.5\\
2426	0.5\\
2427	0.5\\
2428	0.5\\
2429	0.5\\
2430	0.5\\
2431	25.5102040816327\\
2432	0.5\\
2433	0.5\\
2434	0.5\\
2435	0.5\\
2436	0.5\\
2437	0.5\\
2438	0.5\\
2439	0.5\\
2440	0.5\\
2441	0.5\\
2442	0.5\\
2443	0.5\\
2444	0.5\\
2445	0.5\\
2446	0.5\\
2447	0.5\\
2448	0.5\\
2449	0.5\\
2450	0.5\\
2451	25.5102040816327\\
2452	0.5\\
2453	0.5\\
2454	0.5\\
2455	0.5\\
2456	0.5\\
2457	0.5\\
2458	0.5\\
2459	0.5\\
2460	0.5\\
2461	0.5\\
2462	0.5\\
2463	0.5\\
2464	0.5\\
2465	0.5\\
2466	0.5\\
2467	0.5\\
2468	35.7142857142857\\
2469	0.5\\
2470	0.5\\
2471	0.5\\
2472	0.5\\
2473	0.5\\
2474	0.5\\
2475	0.5\\
2476	0.5\\
2477	0.5\\
2478	0.5\\
2479	0.5\\
2480	0.5\\
2481	0.5\\
2482	0.5\\
2483	0.5\\
2484	0.5\\
2485	0.5\\
2486	0.5\\
2487	0.5\\
2488	0.5\\
2489	25.5102040816327\\
2490	0.5\\
2491	0.5\\
2492	0.5\\
2493	0.5\\
2494	0.5\\
2495	0.5\\
2496	0.5\\
2497	0.5\\
2498	0.5\\
2499	0.5\\
2500	0.5\\
2501	0.5\\
2502	0.5\\
2503	0.5\\
2504	0.5\\
2505	0.5\\
2506	0.5\\
2507	0.5\\
2508	35.7142857142857\\
2509	0.5\\
2510	0.5\\
2511	0.5\\
2512	0.5\\
2513	0.5\\
2514	0.5\\
2515	0.5\\
2516	0.5\\
2517	0.5\\
2518	0.5\\
2519	0.5\\
2520	0.5\\
2521	0.5\\
2522	0.5\\
2523	0.5\\
2524	0.5\\
2525	0.5\\
2526	0.5\\
2527	0.5\\
2528	35.7142857142857\\
2529	0.5\\
2530	0.5\\
2531	0.5\\
2532	0.5\\
2533	0.5\\
2534	0.5\\
2535	0.5\\
2536	0.5\\
2537	0.5\\
2538	0.5\\
2539	0.5\\
2540	0.5\\
2541	0.5\\
2542	0.5\\
2543	0.5\\
2544	0.5\\
2545	0.5\\
2546	0.5\\
2547	0.5\\
2548	0.5\\
2549	0.5\\
2550	25.5102040816327\\
2551	0.5\\
2552	0.5\\
2553	0.5\\
2554	0.5\\
2555	0.5\\
2556	0.5\\
2557	0.5\\
2558	0.5\\
2559	0.5\\
2560	0.5\\
2561	0.5\\
2562	0.5\\
2563	0.5\\
2564	0.5\\
2565	0.5\\
2566	0.5\\
2567	35.7142857142857\\
2568	0.5\\
2569	0.5\\
2570	0.5\\
2571	0.5\\
2572	0.5\\
2573	0.5\\
2574	0.5\\
2575	0.5\\
2576	0.5\\
2577	0.5\\
2578	0.5\\
2579	0.5\\
2580	0.5\\
2581	0.5\\
2582	0.5\\
2583	0.5\\
2584	0.5\\
2585	0.5\\
2586	0.5\\
2587	0.5\\
2588	25.5102040816327\\
2589	0.5\\
2590	0.5\\
2591	0.5\\
2592	0.5\\
2593	0.5\\
2594	0.5\\
2595	0.5\\
2596	0.5\\
2597	0.5\\
2598	0.5\\
2599	0.5\\
2600	0.5\\
2601	0.5\\
2602	0.5\\
2603	0.5\\
2604	0.5\\
2605	0.5\\
2606	0.5\\
2607	35.7142857142857\\
2608	0.5\\
2609	0.5\\
2610	0.5\\
2611	0.5\\
2612	0.5\\
2613	0.5\\
2614	0.5\\
2615	0.5\\
2616	0.5\\
2617	0.5\\
2618	0.5\\
2619	0.5\\
2620	0.5\\
2621	0.5\\
2622	0.5\\
2623	0.5\\
2624	0.5\\
2625	0.5\\
2626	0.5\\
2627	35.7142857142857\\
2628	0.5\\
2629	0.5\\
2630	0.5\\
2631	0.5\\
2632	0.5\\
2633	0.5\\
2634	0.5\\
2635	0.5\\
2636	0.5\\
2637	0.5\\
2638	0.5\\
2639	0.5\\
2640	0.5\\
2641	0.5\\
2642	0.5\\
2643	0.5\\
2644	0.5\\
2645	0.5\\
2646	0.5\\
2647	0.5\\
2648	0.5\\
2649	25.5102040816327\\
2650	0.5\\
2651	0.5\\
2652	0.5\\
2653	0.5\\
2654	0.5\\
2655	0.5\\
2656	0.5\\
2657	0.5\\
2658	0.5\\
2659	0.5\\
2660	0.5\\
2661	0.5\\
2662	0.5\\
2663	0.5\\
2664	0.5\\
2665	0.5\\
2666	35.7142857142857\\
2667	0.5\\
2668	0.5\\
2669	0.5\\
2670	0.5\\
2671	0.5\\
2672	0.5\\
2673	0.5\\
2674	0.5\\
2675	0.5\\
2676	0.5\\
2677	0.5\\
2678	0.5\\
2679	0.5\\
2680	0.5\\
2681	0.5\\
2682	0.5\\
2683	0.5\\
2684	0.5\\
2685	0.5\\
2686	0.5\\
2687	35.7142857142857\\
2688	0.5\\
2689	0.5\\
2690	0.5\\
2691	0.5\\
2692	0.5\\
2693	0.5\\
2694	0.5\\
2695	0.5\\
2696	0.5\\
2697	0.5\\
2698	0.5\\
2699	0.5\\
2700	0.5\\
2701	0.5\\
2702	0.5\\
2703	0.5\\
2704	0.5\\
2705	0.5\\
2706	0.5\\
2707	0.5\\
2708	0.5\\
2709	25.5102040816327\\
2710	0.5\\
2711	0.5\\
2712	0.5\\
2713	0.5\\
2714	0.5\\
2715	0.5\\
2716	0.5\\
2717	0.5\\
2718	0.5\\
2719	0.5\\
2720	0.5\\
2721	0.5\\
2722	0.5\\
2723	0.5\\
2724	0.5\\
2725	25.5102040816327\\
2726	0.5\\
2727	0.5\\
2728	0.5\\
2729	0.5\\
2730	0.5\\
2731	0.5\\
2732	0.5\\
2733	0.5\\
2734	0.5\\
2735	0.5\\
2736	0.5\\
2737	0.5\\
2738	0.5\\
2739	0.5\\
2740	0.5\\
2741	0.5\\
2742	0.5\\
2743	0.5\\
2744	0.5\\
2745	25.5102040816327\\
2746	0.5\\
2747	0.5\\
2748	0.5\\
2749	0.5\\
2750	0.5\\
2751	0.5\\
2752	0.5\\
2753	0.5\\
2754	0.5\\
2755	0.5\\
2756	0.5\\
2757	0.5\\
2758	0.5\\
2759	0.5\\
2760	0.5\\
2761	25.5102040816327\\
2762	0.5\\
2763	0.5\\
2764	0.5\\
2765	0.5\\
2766	0.5\\
2767	0.5\\
2768	0.5\\
2769	0.5\\
2770	0.5\\
2771	0.5\\
2772	0.5\\
2773	0.5\\
2774	0.5\\
2775	0.5\\
2776	0.5\\
2777	0.5\\
2778	0.5\\
2779	0.5\\
2780	0.5\\
2781	25.5102040816327\\
2782	0.5\\
2783	0.5\\
2784	0.5\\
2785	0.5\\
2786	0.5\\
2787	0.5\\
2788	0.5\\
2789	0.5\\
2790	0.5\\
2791	0.5\\
2792	0.5\\
2793	0.5\\
2794	0.5\\
2795	0.5\\
2796	0.5\\
2797	25.5102040816327\\
2798	0.5\\
2799	0.5\\
2800	0.5\\
};
\end{axis}
\end{tikzpicture}%

%% file: ap_step_w_ls.tex
\def\xline{100};
\def\radius{100};

\begin{scope}[scale=0.26]


\clip ({(cos(\thetainit)-0.15)*\radius},{(sin(\thetainit)-0.02)*\radius}) rectangle ({1.05*\radius},{0.02*\radius});

\draw[fill=gray!25!white] (0,0) circle (\radius);
\draw (\xline,-{2*\radius})--(\xline,{2*\radius});
\marker{(\radius,0)}{0.01};
\node[right] at (\radius,0) {$x^\star$};
\edef\xoneinit{\radius*cos(\thetainit)};
\edef\xtwoinit{\radius*sin(\thetainit)};

\xdef\maxlssteps{6};

\node[above] at (0.925*\radius,{(sin(\thetainit)-0.1)*\radius}) {$C$};
\node[above] at (0.99*\radius,{(sin(\thetainit)-0.1)*\radius}) {$D$};

\pgfmathparse{1*\xoneinit}
\edef\xone{\pgfmathresult}
\pgfmathparse{1*\xtwoinit}
\edef\xtwo{\pgfmathresult}

\marker{(\xone,\xtwo)}{0.01};
\node[left] at (\xone,\xtwo) {$x^1$};

\edef\pxone{\xline};
\pgfmathparse{1*\xtwo}
\edef\pxtwo{\pgfmathresult}
\marker{(\pxone,\pxtwo)}{0.01};

\pgfmathparse{sqrt(\pxone*\pxone+\pxtwo*\pxtwo)};
\edef\pxnorm{\pgfmathresult};
\pgfmathparse{\pxone*\radius/\pxnorm};
\xdef\ppxone{\pgfmathresult};
\pgfmathparse{\pxtwo*\radius/\pxnorm};
\xdef\ppxtwo{\pgfmathresult};
\marker{(\ppxone,\ppxtwo)}{0.01};

\draw[-latex] (\xone,\xtwo)--(\pxone,\pxtwo);
\draw[-latex] (\pxone,\pxtwo)--(\ppxone,\ppxtwo);
\draw[-latex] (\xone,\xtwo)--(\ppxone,\ppxtwo);

\pgfmathparse{\ppxone-\xone};
\edef\resone{\pgfmathresult};
\pgfmathparse{\ppxtwo-\xtwo};
\edef\restwo{\pgfmathresult};
\pgfmathparse{sqrt(\resone*\resone+\restwo*\restwo)};
\edef\resnorm{\pgfmathresult};

\pgfmathparse{1/\resnorm};
\edef\step{\pgfmathresult}
\draw[red] (\xone,\xtwo)--({\xone+\step*\radius*\resone},{\xtwo+\step*\radius*\restwo});

\input{line_search.tex}

\end{scope}